\documentclass[12pt]{amsart}
\usepackage{}
\usepackage{amsmath}
\usepackage{amsfonts}
\usepackage{amssymb}
\usepackage[all,cmtip]{xy}           %xypic macro for latex2.09
\usepackage{bbm}
\usepackage{bbding}
\usepackage{txfonts}
\usepackage[shortlabels]{enumitem}
\usepackage{ifpdf}
\ifpdf
\usepackage[colorlinks,final,backref=page,hyperindex]{hyperref}
\else
\usepackage[colorlinks,final,backref=page,hyperindex,hypertex]{hyperref}
\fi
\usepackage{mathrsfs}
\usepackage{amscd}
\usepackage{tikz-cd}
\usepackage[active]{srcltx}
\usepackage{tikz}
\usetikzlibrary{calc}
\usetikzlibrary{arrows,shapes,chains}

\usepackage[colorinlistoftodos]{todonotes}

\makeatletter

\newcommand{\st}{\stackrel}

\newcommand{\CA}{\mathcal{A} }

\newcommand{\CF}{\mathcal{F} }

\newcommand{\CQ}{\mathcal{Q} }

\newcommand{\CT}{\mathcal{T} }
\newcommand{\CX}{\mathcal{X} }
\newcommand{\CY}{\mathcal{Y} }

\newcommand{\EE}{\mathbf{E}}
\newcommand{\D}{\mathbf{D}}

\newcommand{\CI}{\mathrm{inj} }
\newcommand{\CP}{\mathrm{proj} }

\newcommand{\PP}{{\mathrm{Proj}}}
\newcommand{\II}{{\mathrm{Inj}}}

\newcommand{\LMod}{{  \rm{Mod}}}

\newcommand{\lmod}{ {\rm{{mod}}}}

\newcommand{\add}{{\rm{add}}}
\newcommand{\Add}{{\rm{Add}}}

\newcommand{\Ggldim}{{\rm{Ggldim}}}
\newcommand{\gldim}{{\rm{gldim}}}
\newcommand{\domdim}{{\rm{domdim}}}

\newcommand{\pd}{{\rm{pd}}}
\newcommand{\id}{{\rm{id}}}

\newcommand{\Gpd}{{\rm{Gpd}}}
\newcommand{\GP}{{\mathrm{GProj}}}
\newcommand{\Gp}{{\mathrm{Gproj}}}

\newcommand{\Coker}{{\rm{Coker}}}

\newcommand{\op}{{\rm{op}}}

\newcommand{\Hom}{{\rm{Hom}}}
\newcommand{\Ext}{{\rm{Ext}}}
\newcommand{\End}{{\rm{End}}}

\newcommand{\Perpo}{{}^{\underline{\perp}_0}}

\newcommand{\Perpn}{{}^{\underline{\perp}_{n}}}

\theoremstyle{plain}
\newtheorem{thm}{Theorem}[section]
\newtheorem{cor}[thm]{Corollary}
\newtheorem{lem}[thm]{Lemma}

\newtheorem{prop}[thm]{Proposition}

\theoremstyle{definition}
\newtheorem{Def}[thm]{Definition}
\newtheorem{example}[thm]{Example}

\newtheorem{remark}[thm]{Remark}

\newtheorem{fremark}[thm]{Final Remark}

\theoremstyle{plain}

\theoremstyle{definition}

\numberwithin{equation}{section}

\newtheorem{sstheorem}{Theorem}

\begin{document}

\title[Characterizing Higher Auslander-Gorenstein Algebras]{Characterizing higher Auslander(-Gorenstein) Algebras}

\author[Ding, Keshavarz, and Zhou]{Zhenhui Ding, Mohammad Hossein Keshavarz,  and Guodong  Zhou}
% \address{School of Mathematics and Statistics, Nantong University, Seyuan Road, Nantong, Jiangsu, 226019, China}
% \email{keshavarz@ntu.edu.cn}
\address{School of Mathematical Sciences, East China Normal University,  Shanghai 200241, China}
\email{51215500001@stu.ecnu.edu.cn}
\email{keshavarz@ntu.edu.cn}
\email{gdzhou@math.ecnu.edu.cn}
% \address{School of Mathematical Sciences,  Key Laboratory of MEA (Ministry of Education), Shanghai Key Laboratory of PMMP, East China Normal University,  Shanghai 200241, China.}
% \email{keshavarz@ntu.edu.cn}
% \email{gdzhou@math.ecnu.edu.cn}

\renewcommand{\thefootnote}{\alph{footnote}}



\begin{abstract}\label{Abstract}
%It is well known that there is a nice connection between Artin algebras of finite representation type and Auslander algebras. 
It is well known that for Auslander algebras, the category of all (finitely generated) projective modules is an abelian category and this property of  abelianness characterizes Auslander algebras by Tachikawa theorem in 1974. 

Let $n$ be a positive integer. In this paper, by using torsion theoretic methods, we show that $ n $-Auslander algebras can  be characterized by the abelianness of the category of  modules with projective dimension less than $ n $ and two additional properties, extending  the classical Auslander-Tachikawa theorem. By Auslander-Iyama correspondence a categorical characterization of the class of Artin algebras having $ n $-cluster tilting modules is obtained.  

Since higher Auslander algebras are a special case of higher Auslander-Gorenstein algebras, the results are given in the general setting as extending   previous results of Kong.

Higher Auslander-Gorenstein Algebras are also studied from the viewpoint of cotorsion pairs and, as application, we show that they satisfy in two nice equivalences.
\end{abstract}

\maketitle

% \tableofcontents

\section{Introduction}\label{Section: Introduction}
An important problem in the representation theory of algebras is to study  algebras of finite representation type; i.e. Artin algebras that  have only a finite number of isomorphism classes of finitely generated indecomposable modules.
\medskip

In 1971, {Auslander} \cite{Au71} proved a remarkable classical result, called Auslander correspondence, about algebras of finite representation type which introduced a completely new insight to the representation theory of Artin algebras.
He showed that there is a bijective correspondence between the set of Morita equivalence classes of Artin algebras $ \Lambda $ of finite representation type and the set of Morita equivalence classes of Artin algebras $ \Gamma $ such that $ \Gamma $ is an Auslander algebra; i.e. Artin algebras $ \Gamma $ whose global dimensions (denoted by $\gldim (\Gamma)$) are at most two and whose dominant dimensions (denoted by $\domdim(\Gamma)$) are at least two.  See  \cite[Page 52]{Au71} or \cite[Theorem VI.5.7]{ARS}.
\medskip

For an Artin algebra $ \Gamma $, let $ \LMod(\Gamma)$ be the category of all left $ \Gamma $-modules and $ \lmod(\Gamma) $ be the category of all finitely generated left $ \Gamma $-modules. If $ \Gamma $ is an Auslander algebra, then the full subcategory $ \CP (\Gamma)  $ of $ \lmod(\Gamma) $ consisting of all finitely generated projective $ \Gamma $-modules is equivalent to a module category and so is itself  an abelian category \cite[Page 52]{Au71}.
\medskip

This gave  {Tachikawa} the  motivation to study rings $ \Gamma $ such that the full subcategory $  \PP (\Gamma)  $ of $ \LMod(\Gamma) $ consisting of all projective $ \Gamma $-modules is abelian.
He obtained some characterizations of such rings and proved that these algebras are precisely Auslander algebras \cite[Theorem 1]{Ta74}.
For Artin algebras, by combining Auslander's result and Tachikawa's result, we have the following characterization of Auslander algebras.
\medskip

\begin{sstheorem}[\bf Auslander-Tachikawa Theorem]\label{T1} For an Artin algebra $ \Gamma $, the following statements are equivalent:
\begin{itemize}
\item[$ (a) $] $ \Gamma $ is an Auslander algebra;
\item[$ (b) $] the full additive  subcategory $ \CP (\Gamma)  $ of  $ \lmod(\Gamma) $ is an abelian category;
\item[$ (c) $] the full additive   subcategory $ \PP (\Gamma)  $ of $ \LMod(\Gamma) $ is an abelian category.
\end{itemize}
\end{sstheorem}

In 2007,  Iyama in a series of  important papers \cite{Iy07-1, Iy07-2} gave a generalization of Auslander correspondence, called higher Auslander correspondence or Auslander-Iyama correspondence. He  showed that for every positive integer $ n $, there is
a bijective correspondence between the set of equivalence classes of $ n $-cluster tilting subcategories with additive generators (Definition \ref{Def: n-Cluster tilting modules}) of Artin algebras and the set of Morita equivalence classes of  Artin algebras $ \Gamma $ such that $ \gldim (\Gamma) \leq n+1 \leq \rm {domdim}( \Gamma) $, called $ n $-Auslander algebras \cite[Theorem 0.2]{Iy07-2}.
%\vspace{5 mm}
%\medskip

A natural question is whether there exists a characterization  for higher Auslander algebras analogous to the Auslander-Tachikawa Theorem. The first goal of this paper is to answer this question and we have the following characterizations for higher Auslander algebras; see Theorem \ref{Thm: Characterizing Higher Auslander Algebras}.
\begin{sstheorem}[\bf Higher Auslander-Tachikawa Theorem]\label{T2}
Let $ n $ be a positive integer and  $ \Gamma $ be an Artin algebra. Then the following statements are equivalent:
\begin{itemize}
\item[$  (a) $] $ \Gamma $ is an $ n $-Auslander algebra;
 
\item[$  (b) $] $ \CP^{\leq n-1} (\Gamma) $ is an abelian category, $ \CP^{\leq n}(\Gamma)\cap \CI(\Gamma) \subseteq \CP(\Gamma) $, and ${^{ \perp_0} \Gamma } \subseteq  {^{ \perp_{n}} \Gamma} $;

\item[$  (c) $] $ \PP^{\leq n-1} (\Gamma) $ is an abelian category, $ \PP^{\leq n}(\Gamma)\cap \II (\Gamma) \subseteq \PP (\Gamma) $, and ${^{\underline{\perp}_0} \Gamma } \subseteq  {^{\underline{\perp}_{n}} \Gamma} $.
\end{itemize}
\end{sstheorem}
Here, for an Artin  algebra $ \Gamma $, we denote   by $\Gamma^{\op}$ the opposite algebra of $\Gamma$, by  $ \II (\Gamma)  $  the full subcategory of $ \LMod(\Gamma) $ consisting of all injective  $ \Gamma $-modules, by $ \CI(\Gamma) $  the full subcategory of $ \lmod(\Gamma) $ consisting of all finitely generated  injective  $ \Gamma $-modules; for every non-negative integer $ m $, by $ \PP^{\leq m}(\Gamma) $  the full subcategory of $ \LMod(\Gamma) $ consisting of all $ \Gamma $-modules of projective dimension  at most $ m $, and by $ \CP^{\leq m}(\Gamma) $  the full subcategory of $ \lmod(\Gamma) $ consisting of all finitely generated $ \Gamma $-modules of projective dimension at most $ m $. For $m\geq 1$,
we write
$${}^{ \perp_{m}}\Gamma:={}^{ \perp_{[1,m]}}\Gamma =\lbrace M \in  \lmod(\Gamma)  \ \vert \  \Ext^{i}_{\Gamma}(M,\Gamma) = 0, \ \forall  \ 1 \leq i \leq m   \rbrace, $$ and
$${}^{\underline{ \perp}_{m}}\Gamma:={}^{\underline{ \perp}_{[1,m]}}\Gamma = \lbrace M \in  \LMod(\Gamma) \  \vert  \ \Ext^{i}_{\Gamma}(M,\Gamma) = 0, \ \forall  \ 1 \leq i \leq m   \rbrace;$$
we also write    $${}^{ \perp_0}\Gamma := \lbrace M \in  \lmod(\Gamma) \   \vert \  \Hom_{\Gamma}(M,\Gamma) = 0 \rbrace,\
{}^{\underline{ \perp}_0}\Gamma := \lbrace M \in  \LMod(\Gamma) \  \vert \  \Hom_{\Gamma}(M,\Gamma) = 0 \rbrace. $$

Note that only the abelianness property of the category of modules with projective dimension less than $n$ can not characterize $n$-Auslander algebras. For $n\geq 2$, the other conditions are  necessary; see Remarks \ref{Remark: additional conditions  of Auslander-Tachikawa}.

%\medskip

In 1993, Auslander and Solberg \cite{AS93} established a Gorenstein analogue of Auslander correspondence by using the technique of relative homology developed in \cite{AS93a, AS93b, AS93c}. They introduced $\tau$-selfinjective algebras, where $\tau$ is the Auslander-Reiten translation,  as a replacement of representation-finite algebras.  The corresponding analogue of Auslander algebras are characterized as Artin algebras whose  injective dimensions are at most $2$ and whose dominant dimensions are at least $2$, that is, Artin algebras $\Gamma$ such that $  \id({}_{\Gamma}   \Gamma) \leq 2 \leq {\rm {domdim}} (\Gamma)$. Kong  named this class of algebras as quasi-Auslander algebras \cite{Ko14}.  In his work \cite{Ko14} Kong also proved an extension of Auslander-Tachikawa theorem by characterizing quasi-Auslander algebras as Artin algebras over which  the category of finitely generated  Gorenstein projective modules is an abelian category.

To obtain a higher analogue of Auslander-Solberg correspondence, Iyama and Solberg \cite{IS18} introduced pre-cluster tilting subcategories in module categories of Artin algebras which are called $\tau_n$-selfinjective algebras  as a higher version of $\tau$-selfjinjective algebras.  They also gave the notion of $n$-minimal Auslander-Gorenstein algebras  as a higher version of quasi-Auslander algebras, say Artin algebras $\Gamma$ with $  \id({}_{\Gamma}   \Gamma) \leq n+1 \leq {\rm {domdim}} (\Gamma)$.
They showed that there exists a bijective correspondence between Morita equivalence classes of $n$-minimal Auslander-Gorenstein algebras and equivalence classes of finite $n$-precluster tilting subcategories of Artin algebras.
%\medskip

In view of the above Iyama-Solberg correspondence, we could establish a Gorenstein version of Theorem~\ref{T2}  which are also   higher versions of Kong's result \cite[Theorem 2.1]{Ko14}.
\begin{sstheorem}[Theorem \ref{Thm: Characterizing Auslander-Gorenstein algebras I}]\label{T4}
Let $ n $ be a positive integer and  $ \Gamma $ be an Artin  algebra. Then the following statements are equivalent.
\begin{itemize}
\item[$  (a) $] $ \Gamma $ is an $ n $-minimal Auslander-Gorenstein algebra;

\item[$  (b) $] $ \Gp^{\leq n-1}(\Gamma)   $ is an abelian category, $ \CP^{\leq n}(\Gamma)  \cap  \CI(\Gamma)  \subseteq \CP (\Gamma)  $  and ${^{ \perp_0} \Gamma } \subseteq  {^{ \perp_{n}} \Gamma} $;

\item[$  (c) $] $ \GP^{\leq n-1}(\Gamma)  $ is an abelian category, $ \PP^{\leq n}(\Gamma)  \cap  \II (\Gamma)  \subseteq \PP (\Gamma) $  and ${^{\underline{\perp}_0} \Gamma } \subseteq  {}^{\underline{\perp}_{n}} \Gamma  $.
\end{itemize}
\end{sstheorem}
Here for every non-negative integer $ m $, we denote by $\GP^{\leq m}(\Gamma)$  the full subcategory of $\LMod(\Gamma)$ consisting of all $ \Gamma $-modules of Gorenstein projective dimension at most $ m $, by  $\Gp^{\leq m}(\Gamma)$ the full subcategory of $ \lmod(\Gamma) $ consisting of all finitely generated $ \Gamma $-modules of Gorenstein projective dimension at most $m$.

While the results of Section \ref{----Section: Higher Auslander-Gorenstein Algebras} mainly make use of torsion theoretic methods, in Section \ref{----Section: Cotorsion Pairs}, we present some results on higher Auslander-Gorenstein algebras from the viewpoint of cotorsion pairs and relate them to the notion  of torsion-cotorsion triples introduced by Bauer, Botnan, Oppermann and Steen \cite{BBOS}. Especially, by Corollary \ref{Cor: Equivalences for Auslander-Gorenstein Algebras},  if $ n $ is a positive integer and  $ \Gamma $ is an $ n $-minimal Auslander-Gorenstein algebra, then we have that
%we have the following equivalences:
$$ {^{ \perp_0} \Gamma } \simeq \frac{\CI^{\leq 1}(\Gamma)}{\CI^{\leq 1}(\Gamma) \bigcap \CP^{\leq n}(\Gamma)}  \simeq \frac{\CI^{\leq 1}(\Gamma)}{\CI(\Gamma)}.$$

%\medskip

 The layout of this paper is as follows.
 
The first section contains some preliminaries including elementary terminology, basic notions and facts about dominant dimension,  torsion pairs, projectivization and injectivization,  injective resolutions,  and  Gorenstein projective modules.

 Since Theorem~\ref{T2} is a special case of Theorem~\ref{T4},  we study higher Auslander-Gorenstein algebras and prove our main result Theorem~\ref{T4} in the second section. %Obviously, Theorem~\ref{T2} is a special case of this result.}

 In the third section, as an application of our work we study $ \tau_n $-selfinjective algebras and give some characterizations of them in the categorical sense. 
  
The fourth section of this paper, Section \ref{----Section: Cotorsion Pairs}, is devoted to study the relationship between higher Auslander-Gorenstein algebras and cotorsion pairs and in it we present some results on higher Auslander-Gorenstein algebras from the viewpoint of cotorsion pairs and relate them to the notion  of torsion-cotorsion triples. See Corollary \ref{Cor: Equivalences for Auslander-Gorenstein Algebras}.

We end the paper in Section \ref{----Section: Higher Auslander Algebras}, where we restrict our attention to higher Auslander algebras and rewrite some results of previous sections for them. The end example gives a better perspective about some results of the paper.

%\medskip

\section{Preliminaries}\label{----Section: Preliminaries}
In this section, for the convenience of the reader, we collect some definitions and results that will be used throughout the paper.

\subsection{Conventions}\label{Subsection: Conventions}\  

Throughout this article all Artin algebras are assumed to be finitely generated modules over a fixed commutative Artinian ring $R$ with   Jacobson radical $J_{R}$.

Let $ \Gamma $ be an Artin  algebra defined over  $R$. We denote by $ \D $ the canonical  duality $ \Hom_{R}(-,\EE $ $(R/J_R)) : \lmod(\Gamma) \longrightarrow \lmod(\Gamma^{\rm \op} ) $, where $\EE(R/J_R)$ is the injective envelope of $R/J_R$ as an $R$-module.

For a left $ \Gamma $-module $ X $,  $ \EE (X) $ denotes  the injective envelope of $ X $;  $ \pd({}_{\Gamma}X) $ (resp. $ \id({}_{\Gamma}X) $)  is the projective (resp. injective)  dimension of $ X $;  while  $ \Add  (X) $ is the full subcategory of $ \LMod(\Gamma) $ whose objects are direct summands of direct sums of copies of $ X $ and $ \add  (X) $ is  the full subcategory of $ \lmod(\Gamma) $ whose objects are direct summands  of finite direct sums of copies of $ X $. 

Given  a subcategory $\mathcal{C}$ of $\LMod(\Gamma)$,  for each $i\geq 1$, denote  $\mathrm{sub}^i(\mathcal{C})$
 to be the subcategory consisting of
  all modules $X$ which admit   copresentations 
  $$0\to X\to C^0\to \cdots \to C^{i-1}$$
  with $C^0, \dots, C^{i-1}\in \mathcal{C}$, while    $\mathrm{fac}_i(\mathcal{C}) $   is formed by all
  modules $Y$ which admit   presentations 
  $$  C_{i-1}\to \cdots \to C_{0}\to Y\to 0$$
  with $C_0, \dots, C_{i-1}\in \mathcal{C}.$   Also,  by convention, $\mathrm{Sub}^0(\mathcal{C}) =\mathrm{Fac}_0(\mathcal{C}) :=\LMod(\Gamma)$.

  For a module $M\in \lmod(\Gamma)$ and $i\geq 0$, write
  $\mathrm{sub}^i(M):=\mathrm{sub}^i(\mathrm{add}(M))$,
 $\mathrm{fac}_i(M):=\mathrm{fac}_i(\mathrm{add}(M))$,
 $\mathrm{Sub}^i(M):=\mathrm{sub}^i(\mathrm{Add}(M))$,
 $\mathrm{Fac}_i(M):=\mathrm{fac}_i(\mathrm{Add}(M))$.
 In particular,
$\mathrm{sub}(M):=\mathrm{sub}^1(M)$,
 $\mathrm{fac}(M) $ $ :=\mathrm{fac}_1(M)$,
 $\mathrm{Sub}(M):=\mathrm{Sub}^1(M)$,
 $\mathrm{Fac}(M):=\mathrm{Fac}_1(M)$.

\subsection{Dominant Dimension}\label{Subsection: Dominant Dimension}\ 

Let $ n $ be a positive integer and  $ \Gamma $ be a ring.
Recall that for a left $ \Gamma $-module $ M $ if in its  minimal injective resolution
$$ 0 \longrightarrow {_{\Gamma} M} \longrightarrow I^{0} \longrightarrow I^{1} \longrightarrow \cdots, $$
 the first $n$ terms $I^0, \dots, I^{n-1}$ are projective, then   the dominant dimension of $ M $ is at least $ n $, denoted by   $ \mathrm {domdim}   ({_{\Gamma} M}) \geq n$. Hence, the dominant dimension of $ M $ is equal to the smallest $ n $ such that $ I^n $ is not projective, or it is infinite if no such $ n $ exists.

The left dominant dimension of the ring $ \Gamma $ is defined as the dominant dimension of the left regular module $ {_{\Gamma} \Gamma} $. Note that the dominant dimension of a right module and the right dominant dimension of a ring are defined similarly. It is well-known that $ \domdim ( {_{\Gamma} \Gamma})  =  \domdim (\Gamma_{\Gamma}) $ ; see \cite[Theorem 4]{Mu68}. So for the rest of the paper, we will denote both left and right dominant dimension of $ \Gamma $, by $ \domdim (\Gamma) $ and call it the dominant dimension of $ \Gamma $.

The following simple observation will be useful in Section~\ref{----Section: Higher Auslander-Gorenstein Algebras}.
\begin{lem}\label{Lem: dominant implies inclusion} For a ring $\Gamma$, if $ \domdim (\Gamma)\geq n+1, $  then ${^{ \perp_0} \Gamma } \subseteq  {^{ \perp_{n}} \Gamma} $ and ${^{\underline{\perp}_0} \Gamma } \subseteq  {}^{\underline{\perp}_{n}} \Gamma  $.

\end{lem}
\begin{proof}
Let
\begin{equation}\label{Eq:minimal injective resolution} 0 \longrightarrow {{}_\Gamma\Gamma} \longrightarrow I^{0} \longrightarrow \cdots  \longrightarrow I^{n} \longrightarrow I^{n+1}\longrightarrow \cdots\end{equation}
 be the minimal injective resolution of $ {}_\Gamma\Gamma $, then $ I^{j} $ is a projective $ \Gamma $-module for every $ 0 \leq j \leq n$.

To prove that ${^{ \perp_0} \Gamma } \subseteq  {^{ \perp_{n}} \Gamma} $, let $ M \in {^{ \perp_0} \Gamma } $. In the minimal injective resolution \eqref{Eq:minimal injective resolution}, as $I^j$ is projective for every $ 0 \leq j \leq n $, we have   $\Hom_{\Gamma}(M, I^{j})=0$.
Applying the functor $ \Hom_{\Gamma}(M,-) $ to \eqref{Eq:minimal injective resolution} shows that $  \Ext^{i}_{\Gamma}(M,\Gamma)=0 $, for every $ 1 \leq i \leq n $. So $ M \in {^{ \perp_{n}} \Gamma}$.

 The inclusion ${^{\underline{\perp}_0} \Gamma } \subseteq  {}^{\underline{\perp}_{n}} \Gamma  $ can be proved exactly in the same way. 
\end{proof}

\subsection{Torsion Pairs }\label{Subsection: Torsion Pairs}\

Let $ \Gamma  $ be an Artin algebra and $ \mathcal{A}:=\LMod(\Gamma) $. If $ X,Y \in \CA $, then we denote the set of morphisms from $ X $ to $ Y $ in $ \CA $ by $ \CA(X,Y) $.
% or $ \Hom_{\CA}(X,Y) $. 
% All subcategories considered are assumed to be full and closed under isomorphisms and  finite direct sums.  
A torsion pair in $ \CA $ is a pair $ (\CX,\CY) $ of  full subcategories of $\CA$ such that $ \CA (X,Y)=0 $ for all $ X \in \CX $ and $ Y \in \CY $, and these two classes are maximal for this property \cite[Sections 1 and 2]{Di66}, that is,
$$ \CX = {^{\perp_0}\CY}:= \lbrace X \in \CA\  \vert\  \CA (X,Y)=0, \forall \ Y \in \CY \rbrace,$$
$$ \CY = {\CX^{\perp_0}}:= \lbrace Y \in \CA\  \vert\  \CA (X,Y)=0, \forall \ X \in \CX \rbrace.$$
If $ (\CX,\CY) $ is a torsion pair in $ \CA $, $ \CX $ is called a torsion class and $ \CY $ is called a torsion-free class.

A subcategory of $ \CA $ is the torsion class (resp. the torsion-free class) of some torsion pair if and only if it is closed under quotients, direct sums  and extensions (resp.  subobjects, direct products and extensions)  \cite[Theorem 2.3]{Di66}.

A torsion pair $ (\CX,\CY) $ in $ \CA $ is called hereditary if $ \CX $ is also hereditary, i.e. $ \CX $ is closed under subobjects (which is equivalent to $ \CY $ being closed under taking injective envelopes \cite[Theorem 2.9]{Di66}). 
For more details about torsion theories in $ \LMod(\Gamma) $; see \cite[Chapter VI]{St75}.

We would like to mention that a subcategory of $ \lmod(\Gamma) $ is the torsion class (resp. the torsion-free class) of some torsion pair if and only if it is closed under quotients and extensions (resp. submodules and extensions). For more details about torsion theories in $ \lmod(\Gamma) $; see \cite[Subsection 1.1]{AIR} or \cite[Subsection 1.2]{Iy04}.

\subsection{(Co)tiliting Modules}\label{Subsection: (co)tiliting modules}\

We recall some definitions and notations about (co)tilting modules; see \cite{Ho82, Sm84}.

\begin{Def}[Tilting and Cotilting Modules] \label{Def: Tilting and Cotilting Modules}
Let $ \Gamma $ be an Artin algebra. A module $T \in \lmod(\Gamma)$ is said to be a tilting module, if the following three properties are satisfied.
\begin{itemize}
\item[(i)] The projective dimension of $ T $ is at most one;
\item[(ii)] $ \Ext_{\Gamma}^{1}(T,T)=0 $;
\item[(iii)] There is an exact sequence $ 0 \longrightarrow \Gamma \longrightarrow  T'  \longrightarrow T''  \longrightarrow 0 $ with $ T' $ and $ T'' $ in $\add(T)$.
%direct summands of direct sums of copies of $ T $.
\end{itemize}

Dually, a module $T \in \lmod(\Gamma)$ will be called a cotilting module if  $ \D (T) $ is a tilting $\Gamma^{\rm op}$-module.
\end{Def}

\begin{Def}[Ext-projectives and Ext-injectives]\label{Def: Ext-projectives}
Let $ \Gamma $ be an Artin algebra and $ \CX$ be a subcategory of $ \lmod(\Gamma) $ closed under extensions. A module $ X \in \CX $ is said to be $ \Ext $-projective (resp. $ \Ext $-injective) if $ \Ext^{1}_{\Gamma}(X,X')=0 $ (resp. $ \Ext^{1}_{\Gamma}(X',X)=0 $) for all $ X' \in \CX $.
\end{Def}

We denote by $ \mathcal{P} (\CX) $ the direct sum of one copy of each of the indecomposable $ \Ext $-projective objects in $ \CX $ up to isomorphism. 

Similarly, we denote by  $ \mathcal{I}(\CX) $ the direct sum of one copy of each of the indecomposable $ \Ext $-injective objects in $ \CX $ up to isomorphism. Also, the annihilator ideal of $ \CX $ is denoted by $ \mathrm{Ann}_{\Gamma}(\CX) $, i.e. $ \mathrm{Ann}_{\Gamma}(\CX) := \{ \gamma \in \Gamma \vert \gamma X =0, \ \forall  \ X \in \CX \}$.

\subsection{Cotorsion Pairs}\label{Subsection: Cotorsion Pairs}\

For a given class $ \mathcal{C} $ of an abelian category $ \CA $, we let
$$ ^{\perp_1}\mathcal{C}:= \{ A \in \CA \mid \Ext^{1}_{\CA}(A,C)=0, \ \forall \ C \in \mathcal{C} \},$$
$$\mathcal{C }^{\perp_1}:= \{ A \in \CA \mid \Ext^{1}_{\CA}(C,A)=0, \ \forall \ C \in \mathcal{C} \}.$$

Let now $ \CX $ and $ \CY $ be full subcategories of the abelian category $ \CA $. Recall that a pair $ (\CX,\CY) $ is called a cotorsion pair (or cotorsion theory) if $ \CX = {}^{\perp_1}\CY$ and $\CY=\CX^{\perp_1} $, see for instance,  \cite[Definition 7.1.2]{EJ11}. In this case, the class $ \CX $ is called a cotorsion class and the class $ \CY $ is called a cotorsion-free class.

The cotorsion pair $ (\CX,\CY) $ is said to be hereditary, if $ \Ext^{i}_{\CA}(X,Y) = 0$, for all $ i \geq 1, X \in \CX$, and $ Y \in \CY $.

The cotorsion pair $ (\CX,\CY) $ is also called complete, if it satisfies the following two conditions:

\begin{itemize}
\item[$ (i) $] For any object $ A $ of $ \CA $, there exists a short exact sequence $ 0 \longrightarrow yA \longrightarrow xA \longrightarrow A \longrightarrow 0$, where $ yA \in \CY $ and $ xA \in \CX $, i.e. $ \CX $ is a special precovering class;
\item[$ (ii) $] For any object $ A $ of $ \CA $, there exists a short exact sequence $ 0 \longrightarrow A \longrightarrow  \tilde{y}A \longrightarrow  \tilde{x}A \longrightarrow 0$, where $  \tilde{y}A \in \CY $ and $ \tilde{x}A \in \CX $, i.e. $ \CY $ is a special preenveloping class.
\end{itemize}

\subsection{Projectivization and Injectivization }\label{Subsection: Projectivization and injectivization}\

Let  $\Lambda$ be an Artin algebra and $M \in \lmod(\Lambda)$. Denote $\Gamma=\mathrm{End}_\Lambda(M)^{\op}$. Then the natural functor
$$\mathrm{Hom}_{\Lambda}(M, -): \lmod(\Lambda)\to \lmod(\Gamma)$$
restricts to an equivalence $\mathrm{add}(M)\simeq \CP(\Gamma)$; see, for instance, \cite[Proposition II.2.1 (c)]{ARS}.
  
\begin{lem}\label{Lem: projectivisation} Let  $P\in \CP(\Lambda)$ be a finitely generated projective module and $\Gamma=\mathrm{End}_\Lambda(P)^{\op}$.  
\begin{itemize}
    \item[(a)] The functor
$$\mathrm{Hom}_{\Lambda}(P, -):\mathrm{fac}_2(P)\to  \lmod(\Gamma)$$
is an equivalence.
\item[(b)]
The equivalence in (a) extends to another equivalence
$$\mathrm{Hom}_{\Lambda}(P, -):\mathrm{Fac}_2(P)\to  \LMod(\Gamma).$$
\end{itemize}
\end{lem}
 
 \begin{proof}
 The first statement is well known; see, for instance,  \cite[Proposition II.2.5]{ARS}. 
 
 For the second statement, since $P$ is finitely generated, $\mathrm{Hom}_{\Lambda}(P, -)$ commutes with direct sums, thus it establishes an equivalence between 
 $\mathrm{Add}(P)$ and $\PP(\Gamma)$.
 Then the proof proceeds as that of \cite[Proposition II.2.5]{ARS}.  \end{proof}

Dually, let  $\Gamma$ be an Artin algebra and $N \in \lmod(\Gamma)$. Denote $\Lambda=\mathrm{End}_\Gamma(N)^{\op}$. Then the natural functor
$$\D\mathrm{Hom}_{\Gamma}(-, N): \lmod(\Gamma)\to \lmod(\Lambda)$$
restricts to an equivalence $\mathrm{add}(N)\simeq \CI(\Lambda)$.

%\medskip

\begin{lem}\label{Lem: injectivisation} Let  $\CQ\in \CI(\Gamma)$ be a finitely generated injective module and $\Lambda=\mathrm{End}_\Gamma(\CQ)^{\op}$.
\begin{itemize}
    \item[(a)] The functor
$$\D\mathrm{Hom}_{\Gamma}(-, \CQ): \mathrm{sub}^2(\CQ)\simeq \lmod(\Lambda)$$
is an equivalence.
\item[(b)]
Let $I=\mathrm{Hom}_{\Gamma^{\op}}(\D\CQ, \Gamma)$, Then the two functors 
$\D\mathrm{Hom}_{\Gamma}(-, \CQ)$ and $  \mathrm{Hom}_\Gamma(I, -) $ from $\lmod(\Gamma)$ to $ \lmod(\Lambda)$
are naturally isomorphic, which  are furthermore isomorphic to $\D\CQ\otimes_\Gamma-$.  

\item[(c)]
The  equivalence in (a) extends to another equivalence
$$\mathrm{Hom}_\Gamma(I, -)\simeq \D\CQ\otimes_\Gamma-:\mathrm{Sub}^2(\CQ)\to  \LMod(\Lambda).$$

\item[(d)]The two categories $\mathrm{sub}^2(\CQ)$ and $\mathrm{Sub}^2(\CQ)$ are abelian categories.
\end{itemize}
 \end{lem}
 
\begin{proof}
(a) This is the first lemma of \cite[Chapter III, Section 4, Page 48]{Au71}.

(b) It is easy to see that $I$ is a finitely generated projective module, as $I\cong \nu^{-1}\CQ$, where $\nu^{-1}$ is the inverse Nakayama functor. 

The two functors $\D\mathrm{Hom}_{\Gamma}(-, \CQ)$ and $  \mathrm{Hom}_\Gamma(I, -) $ are exact, in particular, right exact, so by Watts' theorem (see,  for example, \cite[Corollary 5.47]{Ro09}), they are isomorphic to a tensor functor. It suffices to check that $\D\mathrm{Hom}_{\Gamma}(\Gamma, \CQ)\cong \D\CQ\cong   \mathrm{Hom}_\Gamma(I, \Gamma) $.
In fact, since  $\D\CQ$ is finitely generated projective, its double $\Gamma$-dual $\mathrm{Hom}_\Gamma(\mathrm{Hom}_{\Gamma^{\op}}(\D\CQ, \Gamma), \Gamma)=\mathrm{Hom}_\Gamma(I, \Gamma)$  is isomorphic to itself. 

(c). Since $I$ is finitely generated projective, $$\mathrm{Hom}_\Gamma(I, -)\simeq \D\CQ\otimes_\Gamma-: \mathrm{Sub}^2(\CQ)\to  \LMod(\Lambda)$$  commutes with direct sums, so the equivalence
$\mathrm{add}(\CQ)\simeq \CI(\Lambda)$ extends to 
$\mathrm{Add}(\CQ)\simeq \II(\Lambda)$. As in the proof of  the first lemma of \cite[Chapter III, Section 4]{Au71}, the latter equivalence extends to
another equivalence
$$  \mathrm{Hom}_\Gamma(I, -)\simeq \D\CQ\otimes_\Gamma-:\mathrm{Sub}^2(\CQ)\to  \LMod(\Lambda).$$

(d) follows from (a) and (c). 
\end{proof}

%\subsection{Injective Resolutions }\label{Subsection: Injective Resolutions} \
\subsection{Iwanaga-Gorenstein algebras}\label{Subsection: Iwanaga-Gorenstein algebras} \

Recall that an Artin algebra $\Gamma$ is called Iwanaga-Gorenstein if its injective dimensions both as a left and as a right $ \Gamma $-module are finite, which  are equal by Zaks \cite[Lemma A]{Za69} ( see also \cite[Proposition 9.1.8]{EJ11}).

For Iwanaga-Gorenstein algebras of injective dimension $n+1$, we have the following proposition.

\begin{prop} \label{Prop: Gorenstein implies inclusion} 
Let $\Gamma$ be an Artinian  Iwanaga-Gorenstein algebra of injective dimension at most $n+1$ with the minimal injective resolution
$$ 0 \longrightarrow {_{\Gamma}\Gamma} \longrightarrow I^{0} \longrightarrow \cdots \longrightarrow I^{n}\longrightarrow I^{n+1}\to 0.$$ Then
$ \CP^{\leq n}(\Gamma)  \cap  \CI(\Gamma) \subseteq \add(\oplus_{i=0}^n I^i). $
\end{prop}
\begin{proof}
If $M \in \CP^{\leq n}(\Gamma)   \cap  \CI(\Gamma) $, then  $ M= \bigoplus_{l=1}^r M_l $, where   $ M_{l} $ is an injective indecomposable $ \Gamma $-module and of projective dimension at most $n$.
%By Theorem \ref{Thm: Iwanaga79},  
By \cite[Theorem 2]{Iw79},
each $ M_l $  is a direct summand of $ I^{j} $ with $0\leq j\leq n+1$. But   $ \pd(M_l)\leq n $,  
%so  by Theorem~\ref{Thm: Iwanaga80}, 
so by \cite[Theorem 2]{Iw80},
$ M_l $ is not a direct summand of $ I^{n+1} $. Thus $M_l$ is  a direct summand of $ I^{j} $ with $ 0 \leq j \leq n $.
This shows that $ \CP^{\leq n}(\Gamma)  \cap  \CI(\Gamma) \subseteq \add(\oplus_{i=0}^n I^i). $
\end{proof}

\subsection{Gorenstein Projective Modules}\label{Subsection: Gorenstein Projective Modules}\

Let $\Gamma$ be an Artin algebra.  Recall that a $\Gamma$-module $ M$ is said to be Gorenstein projective if there exists an exact sequence of projective modules
$$ {\bf P} = \cdots \longrightarrow P_1 \longrightarrow P_0 \longrightarrow P^0 \longrightarrow P^1 \longrightarrow \cdots $$
such that $ M \simeq \mathrm{Im} (P_0 \longrightarrow P^0) $ and such that $ \Hom_{\Gamma}(-,Q) $ leaves the sequence $ {\bf P} $ exact whenever $ Q $ is a projective module; see \cite[Definition 10.2.1]{EJ11}.

The class of all  Gorenstein projective modules is denoted by $ \GP(\Gamma) $  and $ \Gp (\Gamma)  $ denotes the subcategory of finitely generated Gorenstein projective modules.

The Gorenstein projective dimension of a $ \Gamma$-module $ M $ is at most $ n \ (n \in \mathbb{N}_0)$ if $ M $ has a Gorenstein projective resolution of length $ n $, written as  $ \Gpd({}_\Gamma M) \leq n $. 
For every non-negative integer $ m $, we denote by   $ \GP^{\leq m}(\Gamma)  $   the full subcategory of $ \LMod(\Gamma) $ consisting of all $ \Gamma $-modules of   Gorenstein projective dimension  at most $ m $ and by  $\Gp^{\leq m}(\Gamma)  $   the full subcategory of $ \lmod(\Gamma) $ consisting of all finitely generated $ \Gamma $-modules  of   Gorenstein projective dimension  at most $ m $. 

Like the classical case, the supremum of  Gorenstein projective dimensions of all left $ \Gamma $-modules is called left Gorenstein global dimension of $ \Gamma $ and is denoted by $ \Ggldim (\Gamma) $. 
Also,  $ \Ggldim (\Gamma) \leq n $ if and only if $ \Gamma $ is an $ n $-Iwanaga-Gorenstein algebra (i.e. $\id ( {}_{\Gamma} \Gamma) \leq n$ and $ \id (\Gamma_\Gamma) \leq n $). In fact, in this case, we have $\Ggldim (\Gamma) = \id ( {}_{\Gamma} \Gamma) = \id (\Gamma_\Gamma) $;  see  \cite[Corollary 12.3.2]{EJ11} or \cite[Corollary 3.2.6]{Ch10}.

\section{Higher Auslander-Gorenstein Algebras}\label{----Section: Higher Auslander-Gorenstein Algebras}

In this section, we will study higher Auslander-Gorenstein algebras introduced by Iyama and Solberg in \cite{IS18} as a generalization of higher Auslander algebras and give a new characterization of them in terms of the abelianness property of the category of modules of Gorenstein projective dimension less than  $n$. Of course, the results can be considered as a higher generalization of Auslander-Tachikawa theorem.

We recall some basic notions and facts about  higher Auslander-Gorenstein algebras.

Let $ n $ be a positive integer. Recall that an Artin  algebra $ \Gamma $ is an $ n $-minimal Auslander-Gorenstein algebra if
$  \id({}_{\Gamma}   \Gamma) \leq n+1 \leq \rm {domdim} (\Gamma)$ \cite[Definition 1.1]{IS18}. 
$ \Gamma $ is also called a higher Auslander-Gorenstein algebra, if   $ \Gamma $ is an $n$-minimal Auslander-Gorenstein algebra for a certain $ n $.

By \cite[Corollary 5.5]{AR94}, these algebras are Iwanaga-Gorenstein algebras.
It is also easily checked that any $n$-minimal Auslander-Gorenstein algebra $ \Gamma $ is either selfinjective or satisfies $  \id({}_{\Gamma}   \Gamma) = n+1 = \domdim (\Gamma)$;  see \cite[Proposition 4.1]{IS18}.
%\medskip

Notice that, as mentioned in Subsection \ref{Subsection: Gorenstein Projective Modules}, for Iwanaga-Gorenstein algebras, there is no difference between Gorenstein global dimension and injective dimension of $ \Gamma $, i.e. $\mathrm{Ggldim}(\Gamma) = \id({}_{\Gamma}   \Gamma)$, and so $ \Gamma $ is an $ n $-minimal Auslander-Gorenstein algebra if and only if $  \Ggldim (\Gamma) \leq n+1 \leq \rm {domdim} (\Gamma)$.

The following observation is the starting point and the key fact of this paper.

{
\begin{prop}\label{Prop: abelianness for higher Auslander-Gorenstein algebra}
Let $ n $ be a positive integer and   $ \Gamma $ be an $ n $-minimal Auslander-Gorenstein algebra.   Then
$ \Gp^{\leq n-1}(\Gamma)  $ is an abelian category. Moreover,  $ \CP^{\leq n}(\Gamma)  \cap  \CI(\Gamma)  \subseteq \CP (\Gamma)  $  and ${^{ \perp_0} \Gamma } \subseteq  {^{ \perp_{n}} \Gamma} $

\end{prop}

\begin{proof}
By Lemma~\ref{Lem: dominant implies inclusion}, ${\rm {domdim} }(\Gamma)\geq n+1$ implies   ${^{ \perp_0} \Gamma } \subseteq  {^{ \perp_{n}} \Gamma} $.

%\medskip
Let
$$ 0 \longrightarrow {_{\Gamma}\Gamma} \longrightarrow I^{0} \longrightarrow \cdots \longrightarrow I^{n}\longrightarrow I^{n+1}\to 0 $$
be the minimal injective resolution of ${}_\Gamma\Gamma$.
By Proposition~\ref{Prop: Gorenstein implies inclusion},   $ \CP^{\leq n}(\Gamma)  \cap  \CI(\Gamma) \subseteq  \add(\oplus_{i=0}^n I^i)$, while the latter is included in $\CP (\Gamma) $, as $ {\rm {domdim}} (\Gamma)\geq n+1$.  This shows that  $ \CP^{\leq n}(\Gamma)  \cap  \CI(\Gamma) \subseteq \CP (\Gamma) $.

%\medskip

Let  $ \CQ $ be the maximal injective summand of $ \Gamma $ and $\Lambda = \End_{\Gamma} (\CQ)^{\op} $.
Since $ \CQ $ is an injective $ \Gamma $-module, by Lemma~\ref{Lem: injectivisation}, $ \D\Hom_{\Gamma}(-,\CQ) $ establishes  an equivalence between
$ \mathrm{sub}^2(\CQ) $ and $ \lmod(\Lambda) $ and so $ \mathrm{sub}^2(\CQ)  $ is an abelian category.
To complete the proof we show that $$ \Gp^{\leq n-1}(\Gamma)  = \mathrm{sub}^2(\CQ).$$

Let $ X \in \mathrm{sub}^2 (\CQ) $. Thus $ X $ fits into a short exact sequence  $$ 0 \longrightarrow X \longrightarrow E^0  \stackrel{f}{\longrightarrow} E^1 $$ such that $ E^0, E^1  \in {\add (\CQ)} $. Since $\CQ$ is also projective,   $ E^0$ and $ E^1$ are also projective.
Let $ X' $ be the cokernel of the morphism $ E^0 \longrightarrow E^1 $.
Since $ \Ggldim (\Gamma) \leq n+1 $,
$ \Gpd({}_{\Gamma} X') \leq n+1 $, so by \cite[Theorem 2.20]{Ho04-1}  $,  \Gpd({}_{\Gamma} X) \leq n-1$, i.e. $X \in  \Gp^{\leq  n-1}(\Gamma)  $. This shows that $ \Gp^{\leq n-1}(\Gamma)  \supseteq   \mathrm{sub}^2(\CQ)$.

Let $ X $ be a finitely generated  Gorenstein projective $ \Gamma $-module. By definition, $ X $ is a submodule of a projective $ \Gamma $-module $ P $. Since $ \id({}_{\Gamma} \Gamma) \leq n+1 \leq \domdim ( \Gamma )$, $ P $ has the minimal injective resolution
$$ 0 \longrightarrow P \longrightarrow J^0 \longrightarrow J^1 \longrightarrow \cdots \longrightarrow J^{n+1} \longrightarrow 0 $$
such that $ J^i \in \CP (\Gamma)   \cap  \CI(\Gamma)  $, $ i=0, \cdots, n$.
Hence   the injective envelope $  \EE(X)$  of $X$  is a direct summand of $J^0$ and itself is  also a projective module.
In the short  exact sequence $$ 0 \longrightarrow X \longrightarrow \EE(X) \longrightarrow K \longrightarrow 0,$$ we have   $ \Gpd({}_{\Gamma}K) \leq 1 $ and by \cite[Theorem 2.10]{Ho04-1},  there exists an exact sequence
$$ 0 \longrightarrow P' \longrightarrow G \longrightarrow K \longrightarrow 0 $$ such that $ G $ is a finitely generated Gorenstein projective module and $ P' $  is a finitely generated projective $ \Gamma $-module.
Since  $ \domdim (\Gamma) \geq n+1   \geq 2$, both $  \EE(G) $ and $ \EE(\EE(P')/P') $ are projective,  
% so by Theorem~\ref{Thm: Miyachi00}, it follows that
so by \cite[Corollary 1.3]{Mi00}, it follows that
$ \EE(K) $ is a projective $ \Gamma $-module as well.  Hence, $ X $ fits into
\begin{equation}\label{Eq: injective resolution}
0 \to X \to \EE(X) \to \EE(K) \end{equation}
with both $\EE(X)$ and $\EE(K)$ being projective-injective, hence $X$ falls into $   \mathrm{sub}^2(\CQ)$.
%\medskip

Now assume that $ X $ is a $ \Gamma $-module such that $ \Gpd({}_{\Gamma}X) \leq n-1$. By \cite[Theorem 2.10]{Ho04-1}, there exists an exact sequence
\begin{equation}\label{Eq: exact sequence}
0 \longrightarrow P_{n-1} \longrightarrow \cdots \longrightarrow P_1 \longrightarrow G \longrightarrow X \longrightarrow 0 \end{equation}
where $ P_{i} $ is a finitely generated  projective module for every $ 1 \leq i \leq n-1 $ and $ G $ is a finitely generated  Gorenstein projective module. Since $ \Gamma $ is an $ n $-minimal Auslander-Gorenstein algebra, $ \Ggldim (\Gamma) \leq n+1 \leq \domdim (\Gamma) $ and so
every projective $ \Gamma $-module $ P_{i} $ has the minimal injective resolution
\begin{equation}\label{Eq: injective resolutions}  0 \longrightarrow P_i \longrightarrow J^{-i,0} \longrightarrow J^{-i,1} \longrightarrow \cdots \longrightarrow J^{-i,n+1} \longrightarrow 0 \end{equation}
such that $ J^{-i,k} \in \CP (\Gamma) \  \cap \  \CI(\Gamma)  $, for every $ k=0, \dots, n$.
Applying \cite[Corollary 1.3]{Mi00}
%Theorem \ref{Thm: Miyachi00}  
to the exact sequence \eqref{Eq: exact sequence} together with \eqref{Eq: injective resolution} and \eqref{Eq: injective resolutions} then immediately gives the required result.

We have shown that $ \Gp^{\leq n-1}(\Gamma)  =  \mathrm{sub}^2(\CQ)$.
\end{proof}
}

In fact, we can  show that the properties obtained in Proposition~\ref{Prop: abelianness for higher Auslander-Gorenstein algebra} are also sufficient to characterize higher Auslander-Gorenstein algebras. We need the following simple lemmas which are of independent interest. They should be known to the expert. However, we could not find them in the literature, so we include a proof.

{
 \begin{lem} \label{Lem: Coincidence of kernels} Let $\Gamma$ be an Artin algebra. If a full additive category $\mathcal{C}$ of $\lmod(\Gamma)$ has kernels and contains $\Gamma$, then the kernels in $\mathcal{C}$ coincides with the kernels in $\lmod(\Gamma)$ and the inclusion functor from $\mathcal{C}$ to  $\lmod(\Gamma)$ is left exact.

 \end{lem}

\begin{proof}Let $ f: X \longrightarrow Y $ be a morphism  in $\mathcal{C}$.
Let  $k': K' \longrightarrow X  $ be the kernel  of $ f $ in  $\mathcal{C}$. We consider the following two complexes and  the natural  isomorphism between them.
\[\xymatrix{ 0 \ar[r] & \Hom_{\Gamma}(\Gamma, K') \ar[r] \ar@{.>}[d]_{\simeq} &  \Hom_{\Gamma}(\Gamma, X) \ar[r] \ar[d]_{\simeq} &  \Hom_{\Gamma}(\Gamma, Y) \ar[d]_{\simeq} \\ 0 \ar[r] & K' \ar[r]^{k'}  & X \ar[r]  & Y }\]
Since   $ \Gamma \in   \mathcal{C}$ the upper complex is an exact sequence   of abelian groups, so the lower complex is also exact in the category of abelian groups. Thus $ k' $ is also the kernel of $ f $ in $ \lmod(\Gamma) $.
\end{proof}

\begin{lem} \label{Lem: bounded Ggldim}
Let $\Gamma$ be an Artin algebra and  $n\geq 1, 0\leq k\leq n+1$. The following statements are equivalent:

\begin{itemize}

\item[$(a)$]$\mathrm{Ggldim}(\Gamma) \leq n+1$;

\item[$(b)$]$\mathrm{sub}^k(\Gp^{\leq n+1-k}(\Gamma))\subseteq \Gp^{\leq n+1-k}(\Gamma)$;

\end{itemize}

Specially, the above statements are equivalent to each of the following:
\begin{itemize}
\item[$(b')$]$\Gp^{\leq n}(\Gamma)$ is closed under taking submodules;

\item[$(b'')$]$\Gp^{\leq n-1}(\Gamma)$ is closed under taking kernels in $\lmod(\Gamma)$.

\end{itemize}

\end{lem}

\begin{proof}
$(a) \Rightarrow (b) $. We prove by induction on $k$. For $k=0$, this is trivial. For $k=1$, let $X\in\Gp^{\leq n}(\Gamma)$ and $Y$ be a submodule of $X$. Form the short exact sequence $0\to Y\to X\to Z\to 0  $.   As $ \Ggldim (\Gamma) \leq n+1 $,
  $\Gpd({}_{\Gamma} Z) \leq n+1$, and combining with $\Gpd({}_{\Gamma} X) \leq n$, we get $\Gpd({}_{\Gamma} Y) \leq n$ by \cite[Theorem 2.20]{Ho04-1}.

  Suppose $(b)$ holds for $k$ ($1\leq k\leq n$). Then for $Y\in \mathrm{sub}^{k+1}(\Gp^{\leq n-k}(\Gamma))$, we have an exact sequence
  $$ 0 \longrightarrow Y \longrightarrow X^{0} \stackrel{d^0}{\longrightarrow} X^{1} \longrightarrow \cdots \longrightarrow X^k $$
  where $X^0,\dots,X^k\in\Gp^{\leq n-k}(\Gamma)$. Let $K^1 :=\Coker(d^0)$. Hence
  $$K^1\in\mathrm{sub}^{k}(\Gp^{\leq n-k}(\Gamma))\subseteq \mathrm{sub}^{k}(\Gp^{\leq n-k+1}(\Gamma))\subseteq\Gp^{\leq n-k+1}(\Gamma)$$
  where the last inclusion holds by inductive assumption.
  Combining with $\Gpd({}_{\Gamma} X^0) \leq n-k$, we get $\Gpd({}_{\Gamma} Y) \leq n-k$ by \cite[Theorem 2.20]{Ho04-1}.

$(b) \Rightarrow (a)$. Let  $ X \in \lmod(\Gamma) $ and $ k  \geq 1 $. Given a resolution of $ X $ by finitely generated projective modules:
  $$ P_{k-1} \st{\beta} \longrightarrow P_{k-2} \longrightarrow\cdots\longrightarrow X \longrightarrow 0, $$
  since $P_0,\dots,P_{k-1}\in \CP(\Gamma) \subseteq \Gp^{\leq n+1-k}(\Gamma)$, we have $K:=\ker(\beta)\in\mathrm{sub}^k(\Gp^{\leq n+1-k}(\Gamma))\subseteq \Gp^{\leq n+1-k}(\Gamma)$. So  $ \Gpd({}_{\Gamma} X ) \leq n+1 $.
   This shows that $ \Ggldim (\Gamma) \leq n+1 $.
\end{proof}

\begin{prop}\label{Prop: Characterization via Small Modules}
Let $ n $ be a positive integer and  $ \Gamma $ be an Artin  algebra. Then $ \Gamma $ is an $ n $-minimal Auslander-Gorenstein algebra  if
$ \Gp^{\leq n-1}(\Gamma)  $ is an abelian category, $ \CP^{\leq n}(\Gamma)  \cap  \CI(\Gamma)  \subseteq \CP (\Gamma)  $, ${^{ \perp_0} \Gamma } \subseteq  {^{ \perp_{n}} \Gamma} $.
\end{prop}

Before giving the proof, we would like to mention that torsion pairs play  an important rule in the proof of  Tachikawa theorem (\cite[Theorem 1]{Ta74}), Iyama's theorem (\cite[Theorem 3.3]{Iy04}), and Kong's theorem (\cite[Theorem 2.1]{Ko14}). We also follow this procedure and prove the result in three steps.

\begin{proof}

{\bf Step I}.
{We show that $ \Ggldim (\Gamma) \leq n+1 $.}

Since  $ \Gp^{\leq n-1} (\Gamma) $ is abelian, it has kernels which, by Lemma~\ref{Lem: Coincidence of kernels}, coincides with the kernels  in $\lmod(\Gamma) $.    Hence, $ \Ggldim (\Gamma) \leq n+1 $ follows from Lemma~\ref{Lem: bounded Ggldim}.

%\medskip

{\bf Step II}.
{We show that $ \Gp^{\leq n} (\Gamma) $ is closed under taking injective envelopes, i.e. $ X \in  \Gp^{\leq n}(\Gamma)$ $ \Rightarrow \EE (X) \in \Gp^{\leq n}(\Gamma) $.}
%\medskip

The idea is to show that $ \Gp^{\leq n}(\Gamma)$ is the torsion-free class of some hereditary torsion pair, and so is closed under taking injective envelopes \cite[Theorem 2.9]{Di66}.

Obviously the subcategory $ \Gp^{\leq n}(\Gamma)  $ is closed under extensions. By Lemma \ref{Lem: bounded Ggldim} ($b'$), it is also closed under taking submodules. Hence, the subcategory $ \Gp^{\leq n}(\Gamma)  $ is closed under taking  extensions and submodules and so is the torsion-free class of some torsion pair $ (\CT, \Gp^{\leq n}(\Gamma) ) $ in $ \lmod(\Gamma) $ \cite[Theorem 2.3]{Di66}.

%\medskip

We now show that $  {^{\perp_0} \Gamma }  =  {^{{\perp}_0} \Gp^{\leq n}(\Gamma)} = \CT$.
Let $ X \in \CT $. Then $ \Hom_{\Gamma}(X, \Gamma)=0$ as $ \Gamma \in \Gp^{\leq n}(\Gamma)$, and so $X\in {^{\perp_0} \Gamma }$. This shows $  {^{\perp_0} \Gamma } \supseteq  \CT$.

On the other hand,  suppose $ X \in {^{\perp_0}} \Gamma $. By  ${}^{\perp_0}   \Gamma  \subseteq  {}^{\perp_{n}}  \Gamma$, $\Ext^{i}_{\Gamma}(X,\Gamma)=0$ for every $ 0 \leq i \leq n $.
Hence, if  $ Y $ is a projective $ \Gamma $-module, then $ \Ext^{i}_{\Gamma}(X,Y)=0$ for every $ 0 \leq i \leq n $. 
If  $ Y $  is a $ \Gamma $-module with $ \pd ({}_\Gamma Y)=1 $, then $ Y $ has a projective resolution as the following.
$$ 0 \longrightarrow P_1 \longrightarrow P_0 \longrightarrow Y \longrightarrow 0 $$
Applying $ \Hom_{\Gamma}(X,-) $ to this exact sequence gives the long exact sequence

\begin{center}
%\vspace{2 mm}

\begin{tikzpicture}[descr/.style={fill=white,inner sep=1.5pt}]
        \matrix (m) [
            matrix of math nodes,
            row sep=1em,
            column sep=2.5em,
            text height=1.5ex, text depth=0.25ex
        ]
        { 0 & \Hom_{\Gamma}(X,P_1) & \Hom_{\Gamma}(X,P_0) & \Hom_{\Gamma}(X,Y) \\
            & \Ext^{1}_{\Gamma}(X,P_1) & \Ext^{1}_{\Gamma}(X,P_0) & \Ext^{1}_{\Gamma}(X,Y) \\
            & \mbox{}         &                 & \mbox{}         \\
             & \Ext^{n-1}_{\Gamma}(X,P_1) & \Ext^{n-1}_{\Gamma}(X,P_0) & \Ext^{n-1}_{\Gamma}(X,Y) \\
            & \Ext^{n}_{\Gamma}(X,P_1) & \Ext^{n}_{\Gamma}(X,P_0) & \Ext^{n}_{\Gamma}(X,Y) \\            
        };

        \path[overlay,->, font=\scriptsize,>=latex]
        (m-1-1) edge (m-1-2)
        (m-1-2) edge (m-1-3)
        (m-1-3) edge (m-1-4)
        (m-1-4) edge[out=355,in=175,black]   (m-2-2)
        (m-2-2) edge (m-2-3)
        (m-2-3) edge (m-2-4)
        (m-2-4) edge[out=355,in=175,dashed, black]   (m-4-2)
        (m-4-2) edge (m-4-3)
        (m-4-3) edge (m-4-4)
        (m-4-4) edge[out=355,in=175,black] (m-5-2)
        (m-5-2) edge (m-5-3)
        (m-5-3) edge (m-5-4);
\end{tikzpicture}
%\vspace{2mm}
\end{center}
and so $ \Ext^{i}_{\Gamma}(X,Y)=0$ for every $ 0 \leq i \leq n-1 $. 
By applying this argument for finite steps,  we deduce that $ \Hom_{\Gamma}(X,Y) = 0 = \Ext^{1}_{\Gamma}(X,Y)$ for each $ Y \in \CP^{\leq n-1}(\Gamma) $. We also deduce that $ \Hom_{\Gamma}(X,Y)=0 $ for each $ Y \in \CP^{\leq n}(\Gamma) $.

Let now $ Y $ be a $ \Gamma $-module with $ \Gpd({}_{\Gamma}Y) \leq n $. Then, by  \cite[Theorem 2.10]{Ho04-1}, there exists an exact sequence
$$ 0 \longrightarrow K \longrightarrow G \longrightarrow Y  \longrightarrow 0 $$
such that $ G $ is a Gorenstein projective module and $ K $  is a $ \Gamma $-module  of projective dimension less than $ n $. By applying $ \Hom_{\Gamma}(X,-) $ to this exact sequence, we deduce that $ \Hom_{\Gamma}(X,Y)=0 $ for each $ Y \in \Gp^{\leq n}(\Gamma)$ as $ G $ is a submodule of a projective $ \Gamma $-module and the projective dimension of $K$ is less than $n$. We have shown that $  {^{\perp_0} \Gamma } = \CT$.
%\medskip

We now prove that  ${^{\perp_0}} \Gamma$ is closed under taking  submodules.
For every short exact sequence $$ 0 \longrightarrow X' \longrightarrow X \longrightarrow X'' \longrightarrow 0 $$ in $ \lmod(\Gamma) $,  if $X \in {^{\perp_0}} \Gamma $, then $ X'' \in  {^{\perp_0}} \Gamma $.
By the exact sequence $$ 0 \longrightarrow  \Hom_{\Gamma}(X'',\Gamma) \longrightarrow  \Hom_{\Gamma}(X,\Gamma) \longrightarrow  \Hom_{\Gamma}(X',\Gamma) \longrightarrow \Ext^{1}_{\Gamma}(X'',\Gamma)$$ and the fact that  $\Ext^{1}_{\Gamma}(X'',\Gamma)=0$  (as  ${}^{\perp_0}   \Gamma  \subseteq  {}^{\perp_{n}}  \Gamma$), $\Hom_{\Gamma}(X',\Gamma)=0$.
Hence, $ {^{\perp_0}} \Gamma $ is also closed under taking  submodules  and so the pair $({^{\perp_0} \Gamma }, \Gp^{\leq n}(\Gamma))$ is a hereditary torsion pair in $\lmod(\Gamma)$.

This finishes the proof of {\bf Step II}.

%\medskip

{\bf Step III}.
{We show that   $ \domdim (\Gamma) \geq n+1 $.}

By {\bf Step II}, for every $ X \in \Gp^{\leq n} (\Gamma) $, $  \EE(X) \in  \Gp^{\leq n} (\Gamma) $.
On the other hand, by \cite[Theorem 2.2]{Ho04-2}, the projective dimension of $ \EE(X) $ equals to its Gorenstein projective dimension and so  $  \EE(X) \in \CP^{\leq n}(\Gamma)  $. 

Hence, by the assumption, $ \EE(X) $ is a projective $ \Gamma $-module.
Since $ \Gamma  \in \Gp^{\leq n}(\Gamma)$, $ I^0 := \EE(\Gamma) \in \CP (\Gamma) $ and so there is an exact sequence
$$ 0 \longrightarrow \Gamma \longrightarrow I^0 \longrightarrow K_0 \longrightarrow 0 $$
with $ I^0 \in \CP (\Gamma)  $ and $ K_0 \in \CP^{\leq 1}(\Gamma) $.
Hence, $ K_0 \in \Gp^{\leq n}(\Gamma) $ and so $ I^1 := \EE(K_0) \in \CP (\Gamma)  $.
Now by applying this argument for finite steps, we deduce that there is an exact sequence
$$ 0 \longrightarrow {}_\Gamma\Gamma \longrightarrow I^0 \longrightarrow I^1 \longrightarrow \cdots \longrightarrow I^n \longrightarrow I^{n+1} \longrightarrow 0 $$
with $ I^{i} \in \CP (\Gamma)   \cap  \CI (\Gamma) $ for each $ 0 \leq i \leq n $ and hence $ \domdim \ \Gamma \geq n+1 $.
\end{proof}

Propositions \ref{Prop: abelianness for higher Auslander-Gorenstein algebra} and  \ref{Prop: Characterization via Small Modules} could be enhanced to infinitely generated modules,  that is,  we could show the following result.  Since its proof is analogous to that of Propositions \ref{Prop: abelianness for higher Auslander-Gorenstein algebra} and  \ref{Prop: Characterization via Small Modules}, we only give a sketch of proof.

\begin{prop}\label{Prop: Characterization via Big Modules}
Let $ n $ be a positive integer and  $ \Gamma $ be an Artin  algebra. Then  $ \Gamma $ is an $ n $-minimal Auslander-Gorenstein algebra if and only if
$ \GP^{\leq n-1}(\Gamma)  $ is an abelian category, $ \PP^{\leq n}(\Gamma)  \cap   \II (\Gamma)  \subseteq \PP (\Gamma)  $ and $ {\Perpo}\Gamma   \subseteq  {\Perpn}\Gamma   $.
\end{prop}

\begin{proof}
``$  \Rightarrow $".    Let $ \Gamma $ be an  $ n $-minimal Auslander-Gorenstein algebra.  Let $X$ be an arbitrary  injective $ \Gamma $-module. 
%By  Facts~\ref{Facts: big modules} $(b)$,
By \cite[Theorem 3.2]{KMT},
$X$ is a direct sum of indecomposable injective  $ \Gamma $-modules, the latter being finitely generated over an Artin algebra. From this fact, we deduce that each   injective $ \Gamma $-module  is a direct sum of finitely generated injective  $ \Gamma $-modules.

If $ X\in  \PP^{\leq n}(\Gamma)  \cap  \II (\Gamma)  $,  then  by the previous paragraph, each indecomposable    direct  summand of $ X $ lies in $  \CP^{\leq n}(\Gamma)   \cap  \CI (\Gamma) $, which by Proposition~\ref{Prop: abelianness for higher Auslander-Gorenstein algebra} also  falls into $ \CP (\Gamma) $. Therefore, $ \PP^{\leq n}(\Gamma)  \cap  \II (\Gamma)  \subseteq \PP (\Gamma) $.

%\medskip

By a similar  argument used to prove Lemma~\ref{Lem: dominant implies inclusion}, one can also easily show that ${^{\underline{\perp}_0} \Gamma } \subseteq  {\Perpn} \Gamma $.

%\medskip

Let now $ \CQ $ be the maximal injective summand of $ \Gamma $.
The statement that  $ \PP^{\leq n-1}(\Gamma)  =  \mathrm{Sub}^2(\CQ) $ can be proved similarly as in the proof of Proposition~\ref{Prop: abelianness for higher Auslander-Gorenstein algebra}. By Lemma~\ref{Lem: injectivisation} (d), $ \GP^{\leq n-1}(\Gamma) $ is abelian.
%\medskip

``$\Leftarrow$". For this direction, the proof is nearly the same as the proof of the implication ``$\Leftarrow$" in Proposition \ref{Prop: Characterization via Small Modules}  except that in {\bf Step II}, one needs to observe that the subcategory $ \GP^{\leq n}(\Gamma)  $ is also closed under taking direct products and to use  \cite[Theorem 2.3]{Di66} to conclude.
%to replace \cite[Chapter VI, 1.4 Proposition]{ARS}  by  \cite[Theorem 2.3]{Di66}.
\end{proof}

Now by combining Propositions \ref{Prop: abelianness for higher Auslander-Gorenstein algebra}, \ref{Prop: Characterization via Small Modules}, and \ref{Prop: Characterization via Big Modules} we have the following theorem that can be considered as a higher generalization of Auslander-Tachikawa theorem.

\begin{thm} \label{Thm: Characterizing Auslander-Gorenstein algebras I}
Let $ n $ be a positive integer and  $ \Gamma $ be an Artin  algebra. Then the following statements are equivalent:
\begin{itemize}
\item[$  (a) $] $ \Gamma $ is an $ n $-minimal Auslander-Gorenstein algebra;

\item[$  (b) $] $ \Gp^{\leq n-1}(\Gamma)   $ is an abelian category, $ \CP^{\leq n}(\Gamma)  \cap  \CI(\Gamma)  \subseteq \CP (\Gamma)  $  and ${^{ \perp_0} \Gamma } \subseteq  {^{ \perp_{n}} \Gamma} $;

\item[$  (c) $] $ \GP^{\leq n-1}(\Gamma)  $ is an abelian category, $ \PP^{\leq n}(\Gamma)  \cap  \II (\Gamma)  \subseteq \PP (\Gamma) $  and ${^{\underline{\perp}_0} \Gamma } \subseteq  {}^{\underline{\perp}_{n}} \Gamma  $;
\end{itemize}
\end{thm}

Based on the proof of Theorem \ref{Thm: Characterizing Auslander-Gorenstein algebras I} (in fact, Propositions \ref{Prop: abelianness for higher Auslander-Gorenstein algebra}, \ref{Prop: Characterization via Small Modules}, and \ref{Prop: Characterization via Big Modules}) we have the following results.

\begin{cor}\label{Cor: Torsion Theories I}
Let $ n $ be a positive integer, $ \Gamma $ be an $ n $-minimal Auslander algebra,  and $ \CQ $ be the maximal injective summand of $ \Gamma $. Then the pairs $ ({^{ \perp_0} \Gamma } , \Gp^{\leq n} (\Gamma))$ and $ ({^{ \underline{\perp}_0} \Gamma } , \GP^{\leq n} (\Gamma))$ are hereditary torsion pairs. % cogenerated by $ \Gamma $ and also  by $ \CQ $. 
\end{cor}

\begin{proof}
In {\bf Step II} of the proof of Proposition \ref{Prop: Characterization via Small Modules}, we showed that the pair $ ({^{ \perp_0} \Gamma } , \Gp^{\leq n} (\Gamma))$ is a hereditary torsion pair. A similar argument shows that the pair $ ({^{ \underline{\perp}_0} \Gamma } , \GP^{\leq n} (\Gamma))$ is also a hereditary torsion pair. 
\end{proof}

\begin{cor}\label{Cor: Subcategories closed under injective envelopes I}
Let $ n $ be a positive integer and   $ \Gamma $ be an $ n $-minimal Auslander-Gorenstein algebra. Then the category of all  $ \Gamma $-modules (resp. finitely generated $ \Gamma $-modules) with dominant dimension at least two is precisely the subcategory  $  \GP^{\leq n-1}(\Gamma) $ (resp. $  \Gp^{\leq n-1}(\Gamma)   $). Also, for every non-negative integer $ i $, the full subcategories $  \GP^{\leq i}(\Gamma)   $ and  $ \Gp^{\leq i}(\Gamma)  $ are closed under taking injective envelopes.
\end{cor}

\begin{proof} For the first statement, it suffices to observe that a   $ \Gamma $-module (resp.    finitely generated $ \Gamma $-module) with dominant dimension at least two is precisely a module in $\mathrm{Sub}^2(\CQ) $ (resp. $\mathrm{sub}^2(\CQ) $).

To prove the second statement, note that by the proof of Theorem \ref{Thm: Characterizing Auslander-Gorenstein algebras I}, the full subcategory $ \GP^{\leq n}(\Gamma) $ of $ \LMod(\Gamma) $ is closed under taking injective envelopes. Thus for every $0 \leq i \leq n $, if $ X \in \GP^{\leq i}(\Gamma) \subseteq \GP^{\leq n}(\Gamma) $, then $  \EE (X) \in  \GP^{\leq n}(\Gamma)   $. But $ \GP^{\leq n}(\Gamma)     \cap   \II (\Gamma)  \subseteq \PP (\Gamma)  $, and so $ \EE (X) \in \PP (\Gamma)   \subseteq \GP^{\leq i}(\Gamma)$.
For each $ i \geq n+1 $, since $ \Ggldim (\Gamma) \leq n+1 $, $ \GP^{\leq i}(\Gamma)  = \LMod(\Gamma) $ and so clearly is closed under taking injective envelopes.

Similarly, for every non-negative integer $ i $,  $  \Gp^{\leq i}(\Gamma)    $ is also closed under taking injective envelopes.
\end{proof}

Kong used of the notion of submodule categories and characterized Artin algebras whose Gorenstein projective modules form an abelian category; see \cite[Corollary 2.3]{Ko14}.  In the following, we will prove a similar property for higher Auslander-Gorenstein algebras.

\begin{cor}\label{Cor: Characterization via submodule categories I}
Let $ n $ be a positive integer and $ \Gamma $ be an Artin algebra.
If $ \Gamma $ is an $ n $-minimal Auslander-Gorenstein algebra, then
$ \Gp^{\leq n}(\Gamma) = \rm{sub} (\Gamma)$ and  $\GP^{\leq n}(\Gamma) = \rm{Sub} (\Gamma) $.
\end{cor}

\begin{proof}
We prove the first equality,  the proof of the second being similar.
Let $ \Gamma $ be an $ n $-Auslander algebra. If $ X \in  \Gp^{\leq n} (\Gamma) $, by Corollary \ref{Cor: Subcategories closed under injective envelopes I},  $ \EE(X)\in \CP^{\leq n} (\Gamma)  \cap  \CI(\Gamma)\subseteq \CP(\Gamma)$,  so   $ X \in \rm{sub} (\Gamma) $.
On the other hand, if a $ \Gamma $-module $ X $ is a submodule of a projective $ \Gamma $-module, then $ \Hom_{\Gamma}(T, X) = 0 $ for every $ T \in {^{ \perp_0} \Gamma } $. But, by Corollary \ref{Cor: Torsion Theories I}, the pair $ ({^{ \perp_0} \Gamma } , \Gp^{\leq n} (\Gamma))$ is a torsion pair and so $ X \in \Gp^{\leq n} (\Gamma) $. Therefore, $ \Gp^{\leq n} (\Gamma)= \rm{sub} (\Gamma)$.
\end{proof}

A combination of our results and some classical results about torsion pairs gives us a better perspective about higher Auslander-Gorenstein algebras, see Corollaries \ref{Cor: Smalo Result I} and \ref{Cor: Ext-projective Objects I}.

\begin{cor} \label{Cor: Smalo Result I}
Let $ n $ be a positive integer and $ \Gamma $ be an $ n $-minimal Auslander-Gorenstein algebra.
Then the following equivalent statements hold:

\begin{itemize}
\item[$ (a) $] $ {^{ \perp_0} \Gamma } $ is functorially finite;
\item[$ (b) $] $ \Gp^{\leq n} (\Gamma) $ is functorially finite;
\item[$ (c) $] $ {^{ \perp_0} \Gamma } =\mathrm{fac}(X) $ for some $ X \in \lmod(\Gamma) $;
\item[$ (d) $] $ \Gp^{\leq n} (\Gamma) = \mathrm{sub}(Y) $ for some  $ Y \in \lmod(\Gamma) $;
\item[$ (e) $] $ \mathcal{P}({^{ \perp_0} \Gamma }) $ is a tilting $( \Gamma/\mathrm{Ann}_{\Gamma}({^{ \perp_0} \Gamma }))$-module;
\item[$ (f) $] $  \mathcal{I}(\Gp^{\leq n} (\Gamma)) $ is a cotilting $\Gamma$-module;
\item[$ (g) $] $ {^{ \perp_0} \Gamma } = \mathrm{fac} ( \mathcal{P}({^{ \perp_0} \Gamma })) $;
\item[$ (h) $] $ \Gp^{\leq n} (\Gamma) = \mathrm{sub} ( \mathcal{I}(\Gp^{\leq n} (\Gamma))) $.
\end{itemize}
\end{cor}

\begin{proof}
By Corollary \ref{Cor: Torsion Theories I}, the pair  $ ({^{ \perp_0} \Gamma } , \Gp^{\leq n} (\Gamma))$  is a torsion pair.  By \cite[Proposition 5.8]{AR94},  $\Gp^{\leq n} (\Gamma)$ is a functorially finite subcategory of $ \lmod(\Gamma) $ and so by \cite[Theorem]{Sm84} and \cite[Proposition 1.1]{AIR},
% so by Theorem \ref{Thm: Smalo Theorem} that we recalled it in the Subsection \ref{Subsection: (co)tiliting modules}, 
all the  statements hold and also are  equivalent with each other. Note that $ \Gamma \in \Gp^{\leq n} (\Gamma) $ and so $ \mathrm{Ann}_{\Gamma}(\Gp^{\leq n} (\Gamma))=0 $. But, by \cite[Lemma 0.2]{Sm84}, $ \mathrm{Ann}_{\Gamma}({^{ \perp_0} \Gamma }) \neq 0 $.
\end{proof}

\begin{cor} \label{Cor: Ext-projective Objects I}
Let $ n $ be a positive integer, $ \Gamma $ be an $ n $-Auslander algebra, and $ \tau $ and $ \tau^{-} $ be the Auslander-Reiten translations. Then the following statements hold.
\begin{itemize}
\item[(i)] A $ \Gamma $-module $ X $ is $ \Ext $-projective in $ \Gp^{\leq n} (\Gamma) $ if and only if $ X $ is a projective $ \Gamma $-module;
\item[(ii)] An indecomposable $ \Gamma $-module $ X $ is $ \Ext $-projective in $ {^{ \perp_0} \Gamma }$ if and only if $ \tau(X) \in  \Gp^{\leq n} (\Gamma)  $, i.e. $ \Gpd_{\Gamma}(\tau X) \leq n $;
\item[(iii)] An indecomposable $ \Gamma $-module $ X $ is $ \Ext $-injective in $  \Gp^{\leq n} (\Gamma)$ if and only if $ \tau^{-}(X) \in  {^{ \perp_0} \Gamma }  $.
\end{itemize}
\end{cor}

\begin{proof}
By Corollary \ref{Cor: Torsion Theories I}, the pair  $ ({^{ \perp_0} \Gamma } , \Gp^{\leq n} (\Gamma))$  is a torsion pair and so, by \cite[Lemmas 1,2, and 3]{Ho82}, all the statements hold; see also \cite[Corollaries 3.4 and 3.7]{AS81}.
\end{proof}

\begin{cor} \label{Cor: Auslander-Reiten Sequences}
Let $ n $ be a positive integer, $ \Gamma $ be an $ n $-minimal Auslander-Gorenstein algebra, and $ t $ be the idempotent radical corresponding to the torsion pair $ ({^{ \perp_0} \Gamma } , \Gp^{\leq n} (\Gamma))$. Moreover, suppose $ X \in  {^{ \perp_0} \Gamma }  $ and $ Y \in   \Gp^{\leq n} (\Gamma) $ are indecomposable $ \Gamma $-modules. Then the following statements hold.
\begin{itemize}
\item[(i)] If $ X $ is not $ \Ext $-projective in $ {^{ \perp_0} \Gamma }  $, then $ t(\tau X)  $ is indecomposable and for the Auslander-Reiten sequence $ 0 \longrightarrow \tau X \longrightarrow E \longrightarrow X \longrightarrow o $, the induced sequence $ 0 \longrightarrow t(\tau X) \longrightarrow t(E) \longrightarrow X \longrightarrow o $ is the Auslander-Reiten sequence in $ {^{ \perp_0} \Gamma }  $;  
\item[(ii)] If $ Y $ is not $ \Ext $-injective in $  \Gp^{\leq n} (\Gamma)   $, then $ \tau^- Y / t(\tau^- Y)  $ is indecomposable and for the Auslander-Reiten sequence $ 0 \longrightarrow Y \longrightarrow E \longrightarrow \tau^- Y \longrightarrow o $, the induced sequence $ 0 \longrightarrow Y \longrightarrow E / t(E) \longrightarrow  \tau^- Y / t(\tau^- Y)  \longrightarrow o $ is the Auslander-Reiten sequence in $  \Gp^{\leq n} (\Gamma) $;
\end{itemize} 
\end{cor}

\begin{proof}
By \cite[Lemmas 2 and 3]{Ho82}, all the statements hold; see also \cite[Corollaries 3.4 and 3.7]{AS81}.
\end{proof}
%\medskip

\section{{$\tau_n $}-selfinjective  Algebras} \label{----Section: tau-n-Selfinjective Algebras}

As an application of our results in Section \ref{----Section: Higher Auslander-Gorenstein Algebras}, in this section we study $ \tau_n $-selfinjective algebras and characterize them in a categorical sense. Let us recall their definition.

\begin{Def}[{$ \tau_n $-selfinjective algebras}]\label{Def: n-Precluster Tilting Modules}
Let $ \Lambda $ be an Artin  algebra and $\CX$ be a full subcategory of $\lmod(\Lambda)$. The subcategory $ \CX $ is cogenerating if  for every object $ M \in \lmod(\Lambda) $ there exist an object $ X \in \CX $ and a monomorphism $ M \longrightarrow X $. The concept of generating subcategories is defined dually.
Also $ \CX $ is called  $ n $-rigid  if $ \Ext^{i}_{\Lambda}(\CX,\CX) =0$ for each $ 0 < i < n $.

Let now $ \tau_{n}$ and $ \tau_{n}^{-} $ be the $ n $-Auslander-Reiten translations; see  \cite[Subsection 1.4]{Iy07-1} or \cite[Definition 1.1]{Iy11}.
A generating-cogenerating $ n $-rigid subcategory $ \CX $ of $ \lmod(\Lambda) $ is called $ n $-precluster tilting subcategory if it is functorially finite,
$ \tau_{n}(\CX)\subseteq \CX $, and $ \tau_{n}^{-}(\CX)\subseteq \CX $   \cite[Definition 3.2]{IS18}. Moreover, if $ \CX $ admits an additive generator $ M $, we say that $ \CX $ is a finite $ n $-precluster tilting subcategory and $ M $ is an $ n $-precluster tilting module. By following Iyama-Solberg    \cite[Definition 3.4]{IS18}, Artin algebras which have $ n $-precluster tilting modules are also called $ \tau_{n} $-selfinjective.
\end{Def}

Note that it follows immediately from the definition and \cite[Theorem 2.3]{Iy07-1} that if $ \CX $ is an $ n $-cluster tilting subcategory of $ \lmod(\Lambda) $, then it  is  necessarily an $ n $-precluster tilting subcategory of $ \lmod(\Lambda) $. Hence, $ \tau_{n} $-selfinjective algebras are in fact a generalization of $ n $-representation finite algebras in the sense of Darp$\rm{\ddot{o}}$-Iyama in \cite[Definition 2.2]{DI20}, see Definition \ref{Def: n-Cluster tilting modules}.

Now, by combining Iyama-Solberg correspondence and our results (Theorem \ref{Thm: Characterizing Auslander-Gorenstein algebras I}), we can prove a result that gives us a new characterization of the class of Artin algebras having $ n $-precluster tilting modules, i.e. $ \tau_n $-selfinjective algebras.
We need the following lemma which is of independent interest and follows from Lemmas~\ref{Lem: projectivisation} and \ref{Lem: injectivisation}.

 \begin{lem}\label{Lem: properties of Auslander-Iyama correspondence}
 Let $ n $ be a positive integer and  $ \Lambda $ be an    Artin  algebra. Let  $M$ be a  $\Lambda$-module which is generating and cogenerating. Denote $\Gamma=\mathrm{End}_\Lambda(M)^{\op}.$ Let  $\CQ:=\mathrm{Hom}_\Lambda(M, \D\Lambda)\cong \D M$ and  $I=\mathrm{Hom}_{\Gamma^{\op}}(\D\CQ, \Gamma)$. Then the following statements hold:
 \begin{itemize}

     \item[(a)] The functor
     $ \mathrm{Hom}_\Lambda(M,-):\lmod(\Lambda)\to \lmod(\Gamma)$ restricts to an equivalence
     $\mathrm{add}(M)\simeq \CP(\Gamma)$, which
       further restricts to another equivalence
     $\CI(\Lambda)\simeq \CP(\Gamma) \cap \CI(\Gamma)$, so
     $\mathrm{add}(\CQ) $ $ =\CP(\Gamma)  \cap   \CI(\Gamma)$.

     \item[(b)] The right $\Gamma$-module $M_\Gamma$ is projective;  $I\cong  \mathrm{Hom}_\Gamma(M, \Gamma)\cong \nu^{-1}\CQ$, where $\nu^{-1}$ is the inverse Nakayama functor, so
     $ \Lambda \cong \mathrm{End}_\Gamma(\CQ)^{\op}  \cong  \mathrm{End}_\Gamma(I)^{\op}.$

      \item[(c)] The three functors
     $M\otimes_\Gamma-, \mathrm{Hom}_\Gamma(I, -)$ and $\D\mathrm{Hom}_\Gamma(-, \CQ)$ from $\lmod(\Gamma)$ to      $\lmod(\Lambda)$  are  naturally isomorphic. So there exists an adjoint triple
     $$\xymatrix{ \lmod(\Lambda)  \ar@<-2ex>[rrr]_{\mathrm{Hom}_\Lambda(M, -)}\ar@<3ex>[rrr]^{I\otimes_\Lambda-} & & &  \lmod(\Gamma)  \ar[lll]_{M\otimes_\Gamma-\simeq \mathrm{Hom}_\Gamma(I, -)}.}$$

     \item[(d)] The functor
     $ \mathrm{Hom}_\Lambda(M,-):\lmod(\Lambda)\to \lmod(\Gamma)$ is fully faithful with essential image
     $  \mathrm{sub}^2(\CQ)  $,
    so $\mathrm{sub}^2(\CQ)\simeq\lmod(\Lambda) $ is an abelian category.

       \item[(e)] The functor
     $ \mathrm{Hom}_\Lambda(M,-):\LMod(\Lambda)\to \LMod(\Gamma)$  is fully faithful with essential image
     $\mathrm{Sub}^2(\CQ)$,
    so $\mathrm{Sub}^2(\CQ) \simeq \LMod(\Lambda)$ is an abelian category.
 \end{itemize}
 \end{lem}
 
\begin{proof}
The  statement $(a)$ can be proved as in \cite[Lemma V.5.3]{ARS}.
For   $(b)$, the fact that   $M_\Gamma$ is projective  follows from the duality
 $\mathrm{Hom}_\Lambda(-, M)$ between $\mathrm{add}(M)$ and $\CP(\Gamma^{op})$ and the fact that ${}_\Lambda\Lambda\in \mathrm{add}(M)$.

Other assertions deduce from $(a)$ and the fact that $\nu^{-1}$ establishes an equivalence between (finitely generated) projective modules and injective modules.

 For $(c)$, the two functors $\mathrm{Hom}_\Gamma(I, -)$ and $D\mathrm{Hom}_\Gamma(-, Q)$ from $\lmod(\Gamma)$ to      $\lmod(\Lambda)$
 are right exact (in fact exact), so by Watts' theorem, one need to verify that they coincides with $M\otimes_\Gamma-$
 on the left regular module ${}_\Lambda \Lambda$, which is obvious by  $Q=\mathrm{Hom}_\Lambda(M, D\Lambda)\simeq D(M)$ and   $I\simeq \mathrm{Hom}_\Gamma(M, \Gamma)\simeq \nu^{-1}Q$. The existence of the adjoint triple follows from the usual tensor-Hom adjunction.

By $(a)$, the functor
     $ \mathrm{Hom}_\Lambda(M,-):\lmod(\Lambda)\to \lmod(\Gamma)$ restricts to an equivalence
       $\CI(\Lambda)\simeq  \mathrm{add}(Q)$, with quasi-inverse
       $D\mathrm{Hom}_\Gamma(-, Q)$. 
       
    By the first lemma of \cite[Chapter III, Section 4, Page 48]{Au71}. this equivalence can be extended to an  equivalence
     $\lmod(\Lambda)\simeq \mathrm{sub}^2(\CQ)$,  as $Q$ is injective.  This proves $(d)$.

     For $(e)$, by \cite[p. 17, Exercise 8]{CE56} and \cite[Theorem 3.2]{KMT},
     % by Fact~\ref{Facts: big modules} $(a)(b)$,
  the  equivalence
       $\CI(\Lambda)\simeq  \mathrm{add}(Q)$ can be extended to  $\II(\Lambda)\simeq  \mathrm{Add}(Q)$. Again the latter can be extended to $  \LMod (\Lambda) \simeq   \mathrm{Sub}^2(\CQ)$.
 \end{proof}

\begin{thm}\label{Thm: equivalence for precluster tilting}
Let $ \Lambda $ be an Artin  algebra and $ n $ be a positive integer. Then   the following are equivalent:
\begin{itemize}
\item[$(a)$] $ \Lambda $ is $ \tau_{n} $-selfinjective, i.e. it has  $ n $-precluster tilting modules;
\item[$(b)$] $ \lmod(\Lambda) $ is equivalent to $ \Gp^{\leq n-1}(\Gamma)  $, where  $ \Gamma $ is an Artin algebra such that $ \CP^{ \leq n}(\Gamma)\cap \CI(\Gamma) \subseteq \CP (\Gamma) $ and ${^{ \perp_0} \Gamma } \subseteq  {^{ \perp_{n}} \Gamma} $;
\item[$(c)$]    $ \LMod(\Lambda) $ is equivalent to $ \GP^{\leq n-1}(\Gamma)  $, where $ \Gamma $ is an Artin algebra such that $ \PP^{\leq n}(\Gamma)\cap \II(\Gamma) \subseteq \PP (\Gamma) $ and ${^{\underline{ \perp}_0} \Gamma } \subseteq  {^{ \underline{\perp}_{n}} \Gamma} $.
\end{itemize}
\end{thm}

\begin{proof}
$(a)\Rightarrow (b)$ follows from Lemma~\ref{Lem: properties of Auslander-Iyama correspondence} $(d)$ and Proposition~\ref{Prop: abelianness for higher Auslander-Gorenstein algebra}.

$(b)\Rightarrow (a)$.
By Proposition \ref{Prop: Characterization via Small Modules},  $\Gamma$ is an $n$-minimal Auslander-Gorenstein algebra and so by Iyama-Solberg correspondence \cite[Theorem 4.5]{IS18}, there is an Artin algebra $ \Lambda' $ with $ n $-precluster tilting $ \Lambda' $-module $ M' $ such that $\lmod (\Lambda') \simeq   \Gp^{\leq n-1}(\Gamma) \simeq \lmod(\Lambda) $.   Hence, $\Lambda$ and $\Lambda'$ are also Morita equivalent.

The equivalence $(a)\Leftrightarrow (c)$
can be proved similarly with Lemma~\ref{Lem: properties of Auslander-Iyama correspondence} (d) replaced by Lemma~\ref{Lem: properties of Auslander-Iyama correspondence} (e) and  Proposition \ref{Prop: Characterization via Small Modules} by
Proposition~\ref{Prop: Characterization via Big Modules}.
\end{proof}

\begin{remark}
Lemmas \ref{Lem: injectivisation} and \ref{Lem: properties of Auslander-Iyama correspondence}, allowed us to generalize our results from finitely generated modules to all modules, see Theorems \ref{Thm: Characterizing Auslander-Gorenstein algebras I} and \ref{Thm: equivalence for precluster tilting}. Of course, we are not limited to these lemmas and can follow other ways. 

Let $ \CQ $ be the maximal injective summand of Artin algebra $\Gamma$. Then $ \CQ $ is product-complete.  Recall that a $\Gamma$-module $ \CQ $ is  product-complete, if   $ \Add (\CQ) =  {\rm Prod} (\CQ) $, where  $ \mathrm{Prod} (\CQ) $ is the smallest subcategory of $  \LMod(\Gamma)  $ containing $Q$ and closed under products and direct summands.

By  \cite[Proposition 3.5]{KS98}, there is a (locally Noetherian) Grothendieck category $ \mathcal{A} $ and a functor $ F: \mathcal{A} \longrightarrow  \LMod(\Gamma)  $ such that $ F $ induces an equivalence between  $ \II (\mathcal{A})  $  and  $ \Add (\CQ) $.
By the proof of  \cite[Proposition 3.5]{KS98}, the functor $ F $ is a composition of the right adjoint of a quotient functor with an exact functor and so $ F $ is left exact. Hence, $ F $ induces an equivalence between $  \CA $ and $ \mathrm{Sub}^2( \CQ ) $. On the other hand, by a similar argument like the proof of Proposition \ref{Prop: abelianness for higher Auslander-Gorenstein algebra}, $  \GP^{\leq n-1}(\Gamma) =  \mathrm{Sub}^2( \CQ ) $ and so $  \GP^{\leq n-1}(\Gamma) $ is an abelian category.
\end{remark}

%\medskip

{

\section{Cotorsion Pairs}\label{----Section: Cotorsion Pairs}

As mentioned in the introduction, this section is devoted to study the relation between higher Auslander-Gorenstein and cotorsion pairs and we present some results on higher Auslander-Gorenstein algebras from the viewpoint of cotorsion pairs and relate them to the notion  of torsion-cotorsion triples introduced by 
Bauer,  Botnan, Oppermann and Steen \cite{BBOS}.

Let $ n $ and $ i $ be a non-negative integer. For an $n$-minimal Auslander-Gorenstein algebra $ \Gamma $, some results by Martinez Villa show that the pair
$$ ({\rm{dom}}^{\geq i}(\Gamma), \CI^{\leq i}(\Gamma)) $$
is a hereditary cotorsion pair in $ \lmod (\Gamma) $, where   $ {\rm{dom}}^{\geq i}(\Gamma) $ denotes the full subcategory of $  \lmod (\Gamma) $ consisting of all finitely generated $ \Gamma $-modules having the dominant dimension at least $ i $  and $ \CI^{\leq i}(\Gamma) $ denotes  the full subcategory of $  \lmod (\Gamma) $ consisting of all finitely generated $ \Gamma $-modules with injective dimension at most $ i $, see \cite[Lemma 3, Proposition 5, and the corollary before it]{Ma92}.

A result by Marczinzik shows that $ \Gp^{\leq n+1-i}(\Gamma)  =  {\rm{dom}}^{\geq i}(\Gamma) $ (see, \cite[Theorem 2.1]{Ma18}) and so, by \cite[Theorem 10]{ET01} or Salce's lemma and \cite[Proposition 7]{Ma92}, the pair
$$ ( \Gp^{\leq n+1-i}(\Gamma),  \CI^{\leq i}(\Gamma)) $$
is a hereditary complete cotorsion pair in $ \lmod (\Gamma) $ for every integer $ i=0,1, \cdots , n+1 $.

The following results describe the kernel of the above mentioned cotorsion pairs and relates them to the notion of tilting theory.

\begin{prop}\label{Prop: n-tilting modules} 
Let $ n $ be a positive integer, $0\leq k\leq n+1$,  $ \Gamma $ be an $ n $-minimal Auslander-Gorenstein algebra; and $Q$ be its maximal projective-injective summand. Then $T_k:=\Omega^{-k}(\Gamma)\oplus Q$ is a $k$-tilting module and $\add(T_k)=\Gp^{\leq k}(\Gamma)\cap  \CI^{\leq n+1-k}(\Gamma)$. Also for all $0\leq j\leq n+1-k$,
$$\mathrm{sub}^j(T_k)={\rm{dom}}^{\geq j}(\Gamma) = \Gp^{\leq n+1-j}(\Gamma).$$
Especially, for every $0\leq i \leq n+1$, $$\mathrm{sub}^i(T_{n+1-i})={\rm{dom}}^{\geq i}(\Gamma) = \Gp^{\leq n+1-i}(\Gamma).$$
\end{prop}
\begin{proof}
Let
$$ 0 \longrightarrow {}_\Gamma\Gamma \longrightarrow I^0 \longrightarrow \cdots \longrightarrow I^n \longrightarrow I^{n+1} \longrightarrow 0 $$
be the minimal injective resolution of $\Gamma$. Since $I^0, I^1, \cdots, I^n \in \add (\CQ)$, clearly, $\Gamma$ has an $\add(T_k)$-copresentation of length $k$ and for every  $0\leq k\leq n+1$, $\Omega^{-k}(\Gamma) \in \Gp^{\leq k}(\Gamma)\cap  \CI^{\leq n+1-k}(\Gamma)$. Hence, for every $1\leq l \leq n+1$,
$$\Ext^{l}_{\Gamma}(T_k,T_k)\simeq \Ext^{l}_{\Gamma}(\Omega^{-k}(\Gamma),\Omega^{-k}(\Gamma)) \simeq \Ext^{l}_{\Gamma}(\Omega^{k}(\Omega^{-k}(\Gamma)),\Gamma) = \Ext^{l}_{\Gamma}(\Gamma,\Gamma) =0. $$
Therefore, $T_k$ is a $k$-tilting module; see \cite[Page 448]{HU96} or \cite[Section 2]{AC01} for the definition. 
%\cite[Subsection 2.2]{An13}

We now show that for every $0\leq j\leq n+1-k$, $\mathrm{sub}^j(T_k)={\rm{dom}}^{\geq j}(\Gamma) = \Gp^{\leq n+1-j}(\Gamma).$

For every $0\leq j\leq n+1-k$, $T_k \in \Gp^{\leq n+1-j}(\Gamma)$ and so by Lemma \ref{Lem: bounded Ggldim},
$$ \mathrm{sub}^j(T_k) \subseteq  \mathrm{sub}^j (\Gp^{\leq n+1-j}(\Gamma)) \subseteq \Gp^{\leq n+1-j}(\Gamma).$$

Conversely, $ \Gp^{\leq n+1-j}(\Gamma) = {\rm{dom}}^{\geq j}(\Gamma) = \mathrm{sub}^j(\CQ) \subseteq \mathrm{sub}^j(T_k). $

To complete the proof, we must show that for every $ 0 \leq k \leq n+1$, $\add(T_k)=\Gp^{\leq k}(\Gamma)\cap  \CI^{\leq n+1-k}(\Gamma)$.
Clearly, $\add(T_k) \subseteq \Gp^{\leq k}(\Gamma)\cap  \CI^{\leq n+1-k}(\Gamma)$ and so we only need to prove the converse.

Let $X \in \Gp^{\leq k}(\Gamma)\cap  \CI^{\leq n+1-k}(\Gamma)$.
By \cite[Theorem 2.2]{Ho04-2},  $  \Gp^{\leq k}(\Gamma) \cap \CI^{\leq n+1-k}(\Gamma) = \CP^{\leq k} (\Gamma) \cap \CI^{\leq n+1-k}(\Gamma) $ and so $\pd (_{\Gamma}X) \leq k$.
Hence, by \cite[Lemma 2.2 $(a)$]{AC01}, $^{\perp}(T_k^{\perp})\subseteq \CP^{\leq k}(\Gamma) \subseteq \Gp^{\leq k}(\Gamma)$ and so $$ X \in \CI^{\leq n+1-k}(\Gamma) = (\Gp^{\leq k}(\Gamma))^{\perp} \subseteq (^{\perp}(T_k^{\perp}))^{\perp} = T_k^{\perp}.$$
On the other hand, by \cite[Lemmas 2.3 $(a)$ and 2.4 $(a)$]{AC01}, for every $Y \in T_k^{\perp} $, there is a short exact sequence $$ 0 \longrightarrow K_0 \longrightarrow T'_0 \longrightarrow Y \longrightarrow o, $$
where $T'_0 \in \add (T_k) $ and $K_0 \in T_k^{\perp}$. By applying this argument for finite steps, there is an exact sequence
$$ 0 \longrightarrow K_{k-1} \longrightarrow T'_{k-1} \longrightarrow \cdots \longrightarrow T'_0 \longrightarrow Y \longrightarrow o,$$
where $T'_0 , \cdots, T'_{k-1} \in \add (T_k) $. Clearly $K_{k-1} \in X^{\perp_{[k+1,\infty)}}$ and $T'_0, \cdots, T'_{k-1} \in \add(T_k)\subseteq X^{\perp}$. Hence, a dimension shifting argument shows that $Y \in X^{\perp}$. Therefore, $X \in {^{\perp}(T_k^{\perp})}$.
But by \cite[Lemmas 2.3 $(a)$ and 2.4 $(a)$]{AC01} there is a short exact sequence $$ 0 \longrightarrow K \longrightarrow T \longrightarrow X \longrightarrow o, $$
where $T \in \add (T_k) $ and $K \in T_k^{\perp}$. Now, by applying functor $\Hom_{\Gamma}(X,-)$, we deduce that $X \in \add (T_k)$.
\end{proof}

Among the above mentioned cotorsion pairs, one of them has a very nice characteristic and its cotorsion class is a torsion-free class of a hereditary torsion pair. This fact relates higher Auslander-Gorenstein Algebras to the notion  of torsion-cotorsion triples introduced by Bauer,  Botnan, Oppermann and Steen \cite{BBOS} and gives us a nice description of them. 

We first recall some definitions and results in this direction and then by using them we will establish two equivalences for $n$-minimal Auslander-Gorenstein algebras.

\begin{Def}[Torsion Cotorsion Triple]\label{Def: Torsion Cotorsion Triple}
A torsion cotorsion triple in an abelian category $ \CA $ is a triple of subcategories $ (\CT, \CF, \mathcal{D}) $ such that the pair $ (\CT, \CF) $ is a torsion pair and the pair $ (\CF, \mathcal{D}) $ is a cotorsion pair, see \cite[Page 29 before Theorem 2.33]{BBOS}.
\end{Def}

\begin{thm}[Dual of {\cite[Theorem 2.35]{BBOS}}]\label{Thm: BBOS Equivalence}
Let $\Gamma$ be an Artin algebra and $ (\CT, \CF, \mathcal{D}) $ be a torsion cotorsion triple in  $ \lmod (\Gamma)$. Then the inverse of Auslander–Reiten translation defines the following equivalence.
$$ \CT \simeq \frac{\mathcal{D}}{\CI(\Gamma)}.$$
\end{thm}

\begin{proof}
This theorem is the dual version of \cite[Theorem 2.35]{BBOS} and for the convenience of the reader we give the proof.

By \cite[Proposition IV.1.9]{ARS}, we know that the inverse of Auslander-Reiten translation $\tau^{-}$ defines an equivalence 
\[\begin{tikzcd}
	\overline{\lmod} (\Gamma) &&  \underline{\lmod} (\Gamma).
	\arrow["\tau^{-} =  \mathrm{Tr} \D ", from=1-1, to=1-3]
\end{tikzcd}\]
Clearly, $  \frac{\mathcal{D}}{\CI(\Gamma)} $ is a full subcategory of $  \overline{\lmod} (\Gamma) $ and one can easily see that $ \CT $ is a full subcategory of $ \underline{\lmod} (\Gamma)  $. In fact, the class $ \CF $ contains projectives and so $ \Hom_{\Gamma}(\CT , \CP(\Gamma)) =0 $. Thus there is no zero maps between objects in $ \CT $ factor through projective modules.

Now, to complete the proof, we show that a module $ X $ is in $ \mathcal{D} $ if and only if $ \tau^{-} X $ is in $ \CT $. We observe that
$$ X \in  \mathcal{D} \Longleftrightarrow \Ext^{1}_{\Gamma}(\CF,X)=0 \Longleftrightarrow \Hom_{\Gamma}(\tau^- X, \CF)=0 \Longleftrightarrow  \tau^- X \in \CT,$$
where the middle equivalence follows from Auslander-Reiten formula because the injective dimension of every element of $ \mathcal{D} $ by Bauer-Botnan-Oppermann-Steen Correspondence is at most one, see \cite[Corollary III.6.4]{SY11} and \cite[Theorem 2.34]{BBOS}.
\end{proof}

As an immediate result of these results, we have the following equivalences.

\begin{cor}\label{Cor: Equivalences for Auslander-Gorenstein Algebras}

Let $ n $ be a positive integer and  $ \Gamma $ be an $ n $-minimal Auslander-Gorenstein algebra. Then we have the following equivalences.
$$ {^{ \perp_0} \Gamma } \simeq \frac{\CI^{\leq 1}(\Gamma)}{\CI^{\leq 1}(\Gamma) \bigcap \CP^{\leq n}(\Gamma)}  \simeq \frac{\CI^{\leq 1}(\Gamma)}{\CI(\Gamma)}.$$
\end{cor}

\begin{proof}
Since $ \Gamma $ is an $ n $-minimal Auslander-Gorenstein algebra, the triple $ (^{ \perp_0} \Gamma, \Gp^{\leq n}(\Gamma),  \CI^{\leq 1}(\Gamma)) $ is a torsion cotorsion triple, see Corollary \ref{Cor: Torsion Theories I} and the arguments at the beginning of this section. Hence, by \cite[Theorem 2.33]{BBOS}, we have the first equivalence. That is
$$ {^{ \perp_0} \Gamma } \simeq \frac{\CI^{\leq 1}(\Gamma)}{\CI^{\leq 1}(\Gamma) \bigcap \CP^{\leq n}(\Gamma)}. $$
The second equivalence also follows from Theorem \ref{Thm: BBOS Equivalence}. 
%That is $$ {^{ \perp_0} \Gamma }  \simeq \frac{\CI^{\leq 1}(\Gamma)}{\CI(\Gamma)}.$$
\end{proof}

\section{Higher Auslander Algebras}\label{----Section: Higher Auslander Algebras}

As higher Auslander algebras are special cases of higher Auslander-Gorenstein algebras and are of independent interest, in this short section we restrict our attention to these algebras and rewrite some of results of previous sections for them. 

Let $ n $ be a positive integer. Recall that an Artin  algebra $ \Gamma $ is an $ n $-Auslander algebra if its global dimension is at most $ n+1 $ and its dominant dimension is at least $ n+1 $, that is,  $$  \gldim  (\Gamma) \leq n+1 \leq \mathrm{domdim} (\Gamma).$$ 
$ \Gamma $ is also called a higher Auslander algebra, if there is a positive integer $ n $ such that $ \Gamma $ is an $n$-Auslander algebra. Iyama introduced these algebras as a generalization of  classical Auslander algebras  \cite[Page 3]{Iy11}.

As a result of Theorem \ref{Thm: Characterizing Auslander-Gorenstein algebras I}, we have the following result for higher Auslander algebras.

\begin{thm}[\bf Higher Auslander-Tachikawa Theorem]\label{Thm: Characterizing Higher Auslander Algebras}
Let $ n $ be a positive integer and  $ \Gamma $ be an Artin algebra. Then the following statements are equivalent.
\begin{itemize}
\item[$  (a) $] $ \Gamma $ is an $ n $-Auslander algebra;
 
\item[$  (b) $] $ \CP^{\leq n-1} (\Gamma) $ is an abelian category, $ \CP^{\leq n}(\Gamma)\cap \CI(\Gamma) \subseteq \CP(\Gamma) $, and ${^{ \perp_0} \Gamma } \subseteq  {^{ \perp_{n}} \Gamma} $;
\item[$  (c) $] $ \PP^{\leq n-1} (\Gamma) $ is an abelian category, $ \PP^{\leq n}(\Gamma)\cap \II (\Gamma) \subseteq \PP (\Gamma) $, and ${^{\underline{\perp}_0} \Gamma } \subseteq  {^{\underline{\perp}_{n}} \Gamma} $. \qed
\end{itemize}
\end{thm}

\begin{remark}\label{Remark: additional conditions of Auslander-Tachikawa}
For $n\geq 2$, the conditions $ \CP^{\leq n}(\Gamma)\cap \CI(\Gamma) \subseteq \CP(\Gamma) $ and ${^{ \perp_0} \Gamma } \subseteq  {^{ \perp_{n}} \Gamma} $ in Theorem~\ref{Thm: Characterizing Higher Auslander Algebras}(b) are necessary. 
In fact, take an Artin algebra of global dimension $n-1$ which is not an $n$-Auslander algebra. Then  $ \CP^{\leq n-1} (\Gamma)=\lmod(\Gamma) $ is an abelian category, but by Theorem~\ref{Thm: Characterizing Higher Auslander Algebras}, the conditions that $ \CP^{\leq n}(\Gamma)\cap \CI(\Gamma) \subseteq \CP(\Gamma) $  or ${^{ \perp_0} \Gamma } \subseteq  {^{ \perp_{n}} \Gamma} $ don't hold.
\end{remark}

\begin{remark}\label{Remark: Study the abelianness of Subcategories}
 Assume that  $ \Gamma $ is an $ n $-Auslander algebra such that $ \gldim \ \Gamma = n+1 = \domdim \ \Gamma $. By Theorem \ref{Thm: Characterizing Higher Auslander Algebras}, $ \CP^{\leq n-1}(\Gamma)  $ is an abelian category and $ \CP^{\leq i}(\Gamma)  $ is not an abelian category for each $ i=0, \dots, n-2$. Hence, Theorem \ref{Thm: Characterizing Higher Auslander Algebras}  can be used  to study the abelianness of some subcategories of $ \lmod(\Gamma) $.
\end{remark}

Now, like Theorem \ref{Thm: equivalence for precluster tilting}, by combining Auslander-Iyama correspondence and our results (Theorem \ref{Thm: Characterizing Higher Auslander Algebras}), we have the following result that gives us a new characterization of the class of Artin algebras having $ n $-cluster tilting modules, i.e. $ n $-representation finite algebras. We first recall the definition of $ n $-representation finite algebras.

\begin{Def}[$ n $-representation finite]\label{Def: n-Cluster tilting modules}
Let $ \Lambda $ be an Artin algebra and $\mathcal{M}$ be a full subcategory of $\lmod(\Lambda)$.
Recall that $ \mathcal{M} $ is an {\bf $ n $-cluster-tilting }   subcategory of $ \lmod(\Lambda) $ if it is functorially finite in $  \lmod(\Lambda) $ and
$$\begin{array}{rcl} \mathcal{M} &=& \lbrace X \in  \lmod(\Lambda) \  \ \vert\  \ \Ext^{i}_{\Lambda}(X,\mathcal{M})=0, \ \forall \ 0 < i < n   \rbrace\\
     &=& \lbrace  X \in  \lmod(\Lambda)\   \ \vert \   \ \Ext^{i}_{\Lambda}(\mathcal{M},X)=0,   \ \forall \ 0 < i < n  \rbrace.\end{array} $$
Moreover, if $ \mathcal{M} $ admits an additive generator $ M $, i.e.  $ \mathcal{M} = \add (M) $, we say that $ \mathcal{M} $ is a finite $ n $-cluster tilting subcategory and $ M $ is an $ n $-cluster tilting module.

By following Darp$\rm{\ddot{o}}$-Iyama in \cite[Definition 2.2]{DI20}, Artin algebras which have $ n $-cluster tilting modules are also called $ n $-representation finite. We do not assume $ \gldim \Gamma \leq n $ in contrast with several earlier papers; see for instance \cite[Definition 2.2]{IO11}.
\end{Def}

\begin{thm}\label{Thm: equivalence for cluster tilting}
Let $ \Lambda  $ be an Artin algebra and $ n $ be a positive integer. Then the following are equivalent.
\begin{itemize}
\item[$(a)$] $ \Lambda $ is $ n $-representation finite, i.e. it has an  $ n $-cluster tilting module;

\item[$(b)$] $ \lmod(\Lambda) $ is equivalent to $ \CP^{\leq n-1}(\Gamma)  $, where  $ \Gamma $ is an Artin algebra such that $ \CP^{ \leq n}(\Gamma)\cap \CI(\Gamma) \subseteq \CP (\Gamma) $ and ${^{ \perp_0} \Gamma } \subseteq  {^{ \perp_{n}} \Gamma} $;

\item[$(c)$] $ \LMod(\Lambda) $ is equivalent to $ \PP^{\leq n-1}(\Gamma)  $, where $ \Gamma $ is an Artin algebra such that $ \PP^{\leq n}(\Gamma)\cap \II(\Gamma) \subseteq \PP (\Gamma) $ and ${^{\underline{ \perp}_0} \Gamma } \subseteq  {^{ \underline{\perp}_{n}} \Gamma} $. \qed
\end{itemize}
\end{thm}

The following example gives a better perspective about some results of the paper. In the following, for every vertex $i$ of the given quiver $Q$, the corresponding simple, projective, and injective object will be denoted by $S(i)$, $P(i)$, and $I(i)$, respectively. For more details about quiver representations, see for example \cite{SY11}.

\begin{example}\label{Ex: Cyclic Quiver}
Let $\Gamma=KQ/I$ be an algebra with the quiver $Q$:
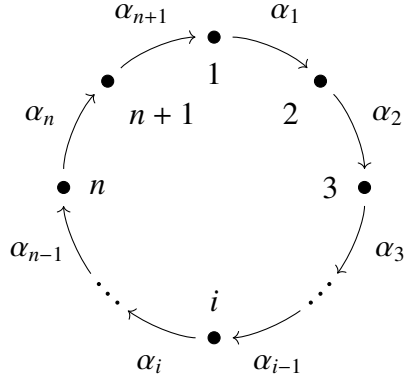
\begin{center}
\begin{tikzpicture}

\foreach \ang\lab\anch in {90/1/north, 45/2/{north east}, 0/3/east, 270/i/south, 180/{n}/west, 135/n+1/{north west}}{
  \draw[fill=black] ($(0,0)+(\ang:2)$) circle (.08);
  \node[anchor=\anch] at ($(0,0)+(\ang:1.8)$) {$\lab$};
 % \draw[->,shorten <=7pt, shorten >=7pt] ($(0,0)+(\ang:3)$).. controls +(\ang+40:1.5) and +(\ang-40:1.5) .. ($(0,0)+(\ang:3)$);
}

% Top part of circle, arrows between different nodes and their labels
\foreach \ang\lab in {90/1,45/2,180/{n},135/n+1}{
  \draw[->,shorten <=7pt, shorten >=7pt] ($(0,0)+(\ang:2)$) arc (\ang:\ang-45:2);
  \node at ($(0,0)+(\ang-22.5:2.5)$) {$\alpha_{\lab}$};
}

% Bottom part of circle, arrows between different nodes and their labels
\draw[->,shorten <=7pt] ($(0,0)+(0:2)$) arc (360:325:2);
\draw[->,shorten >=7pt] ($(0,0)+(305:2)$) arc (305:270:2);
\draw[->,shorten <=7pt] ($(0,0)+(270:2)$) arc (270:235:2);
\draw[->,shorten >=7pt] ($(0,0)+(215:2)$) arc (215:180:2);
\node at ($(0,0)+(0-20:2.5)$) {$\alpha_3$};
\node at ($(0,0)+(315-25:2.5)$) {$\alpha_{i-1}$};
\node at ($(0,0)+(270-20:2.5)$) {$\alpha_i$};
\node at ($(0,0)+(225-25:2.5)$) {$\alpha_{n-1}$};

% Ellipsis
\foreach \ang in {310,315,320,220,225,230}{
  \draw[fill=black] ($(0,0)+(\ang:2)$) circle (.02);
}

\end{tikzpicture}
\end{center}

and $I$ is the admissible  ideal of $KQ$ generated by paths $\alpha_{i+1}\alpha_{i}$ for every $1 \leq i \leq n$. Then:
\begin{enumerate}[\rm(1)]
\item By \cite[Theorem I.10.5]{SY11}, the Nakayama algebra $\Gamma$ is of finite representation type with $2n+3$ indecomposable $\Gamma$-modules. The following table gives us more information about their homological dimensions.
%\medskip

\begin{center}

\setlength{\arrayrulewidth}{1.1pt}
\begin{tabular}{|c|c|c|}
    \hline
    Ind. Modules & \ \ \ \pd \  \ \  & \ \ \ \id \ \ \  \\
    \hline
    \hline
    $P(1)$ & $0$ & $n+1$ \\
    \hline
    $P(2)=I(3)$ & $0$ & $0$ \\
    \hline
    $\vdots$ & $\vdots$ & $\vdots$ \\
    \hline
    $P(n)=I(n+1)$ & $0$ & $0$ \\
    \hline
    $P(n+1)=I(2)$ & $0$ & $0$ \\
    \hline
    $I(1)$ & $n+1$ & $0$ \\
    \hline
    $S(1)$ & $n+1$ & $n+1$ \\
    \hline
    $S(2)$ & $n$ & $1$ \\
    \hline
    $\vdots$ & $\vdots$ & $\vdots$ \\
    \hline
    $S(j) \ (2 \leq j \leq n+1)$ & $n+2-j$ & $j-1$ \\
    \hline
    $\vdots$ & $\vdots$ & $\vdots$ \\
    \hline
    $S(n+1)$ & $1$ & $n$ \\
    \hline
\end{tabular}

\end{center}
%\medskip
\item The global dimension of $\Gamma$ is $n+1$ and $P(i)=I(i+1)$ for
$2 \leq i\leq n$ and $P(n+1)=I(2)$.
\item The minimal injective resolution of $\Gamma$ is as follows:

$0\rightarrow \Gamma \rightarrow(\bigoplus_{i=2}^{n+1}I(i)) \oplus I(2)\rightarrow I(2) \rightarrow I(n+1)  \rightarrow \cdots\rightarrow I(2) \rightarrow I(1)\rightarrow 0$.

%Hence $\Gamma$ is $(n+1)$-Gorenstein.
\item The dominant dimension of $\Gamma$ is $n+1$ and so $\Gamma$ is an $n$-Auslander algebra.
\item The full subcategory $\CP^{\leq n}(\Gamma)$ of $\lmod(\Gamma)$  is generated by indecomposable  $\Gamma$-modules $P(1), \cdots,$ $ P(n+1),  S(2), \cdots , S(n+1)$.
Hence, by Corollary \ref{Cor: Torsion Theories I}, the torsion class ${^{ \perp_0} \Gamma }$ can have only $S(1)$ and $I(1)$ as the idecomposable objects. But $\Hom_{\Gamma}(I(1),P(n))  \neq 0$ and so ${^{ \perp_0} \Gamma } = \add (S(1))$.
\item By Corollary \ref{Cor: Ext-projective Objects I} $(i)$, $ P(1), \cdots, P(n+1) $ are all indecomposable $\Ext$-projective objects in $\CP^{\leq n} (\Gamma)$, i.e.  $\mathcal{P}(\CP^{\leq n} (\Gamma)) = \bigoplus_{i=1}^{n+1}P(i)$.
\item As mentioned in Section \ref{----Section: Cotorsion Pairs}, the pair $ (\CP^{\leq n}(\Gamma),  \CI^{\leq 1}(\Gamma)) $ is a cotorsion pair and so $\CP^{\leq n}(\Gamma) \cap  \CI^{\leq 1}(\Gamma) $ is the class of all $\Ext$-injective objects in $\CP^{\leq n} (\Gamma)$, i.e.  $\mathcal{I}(\CP^{\leq n} (\Gamma)) = (\bigoplus_{i=2}^{n+1}I(i)) \oplus S(2).$ %Also, by Corollary \ref{Cor: Ext-projective Objects I} $(iii)$, $\tau^{-}(S(2))=S(1)$.
\item By Corollary \ref{Cor: Auslander-Reiten Sequences} $(ii)$, the Auslander-Reiten sequence
$$0 \longrightarrow P(1)  \longrightarrow S(1) \oplus P(n+1) \longrightarrow I(1) \longrightarrow 0$$ in $\lmod (\Gamma)$ induces the following Auslander-Reiten sequence in $\CP^{\leq n}(\Gamma)$:
$$0 \longrightarrow P(1) \longrightarrow P(n+1) \longrightarrow S(n+1) \longrightarrow 0.$$
\item  Since $\Ext^{1}_{\Gamma}(S(1),S(1))=0$, $S(1)$ is both $\Ext$-projective and $\Ext$-injective in $^{ \perp_0} \Gamma $. Hence, $\mathcal{P}(^{ \perp_0} \Gamma) = S(1) = \mathcal{I}(^{ \perp_0} \Gamma)$.
\end{enumerate}
\end{example}

\begin{fremark}
Let $ R $ be a commutative artinian ring and $ n $ be a positive integer. In 2017, Iyama and Jasso extended higher Auslander correspondence from Artin $R$-algebras of finite representation type to dualizing $R$-varieties. By their result, i.e. higher Auslander correspondence for dualizing $ R $-varieties, we know that an Artin $ R $-algebra $ \Gamma $ is an  $ n $-Auslander algebra if and only if the category of all finitely generated projective $ \Gamma $-modules is an $ n $-abelian category \cite[Theorem 1.2]{IJ17}. But by some results due to Beligiannis \cite{Be15}, even for $n$-minimal Auslander-Gorenstein algebras, one can see that the category of all finitely generated Gorenstein projective $ \Gamma $-modules is not an $ n $-abelian category in general. Hence,  based on our results, there is a difference between the abelianness of $\Gp^{\leq n-1}(\Gamma)$ and the $n$-abelianness of $\Gp(\Gamma)$ in general.
\end{fremark}

%\medskip

 \textbf{Acknowledgments:}  
 % The authors   were   supported by This work was supported by the China National Key R $\&$ D Program (No. 2024YFA1013803), the Shanghai Key Laboratory of PMMP (No.22DZ2229014) and a grant from IPM (No. 96180070). 
 It is our great pleasure to thank many reseachers for their insightful suggestions and remarks, which include, but not restricted to,   Zhaoyong Huang, Xiaojin Zhang, Bin Zhou, and Shijie Zhu.

\end{document}